\def\ArxivVersion{1}
\def\UsePrebuiltFigures{1}
\documentclass[10pt,journal,twoside]{IEEEtran}

\ifdefined\pdfminorversion
\fi
\ifdefined\pdfobjcompresslevel
\fi
\ifdefined\UsePrebuiltFigures
\fi
\usepackage[dvipsnames]{xcolor}

\usepackage[british]{babel}
\addto\captionsbritish{%
  \renewcommand{\figurename}{Fig.}%
  \renewcommand{\tablename}{TABLE}%
}
\usepackage{textcomp}
\usepackage{xurl}
\usepackage{cite}
\usepackage{enumitem}

\usepackage{amsmath,amssymb,amsfonts,amsbsy}
\usepackage{mathtools}

\usepackage[only,llbracket,rrbracket]{stmaryrd}
\SetSymbolFont{stmry}{bold}{U}{stmry}{m}{n}

\usepackage{mathrsfs}
\usepackage{bm}
\usepackage{nicefrac}
\usepackage{xfrac}
\usepackage{latexsym,theorem}
\usepackage[thinc]{esdiff}
\usepackage{units}

\usepackage{graphicx}
\usepackage{epsfig}
\ifdefined\UsePrebuiltFigures\else
  \usepackage{psfrag}
\fi
\usepackage{eepic}
\usepackage{array}
\usepackage{colortbl}
\usepackage{multirow}
\usepackage{hhline}
\usepackage{booktabs}
\usepackage{float}
\usepackage{multicol}
\usepackage{balance}
\usepackage[strict]{changepage}
\usepackage{cuted}
\usepackage[rightcaption]{sidecap}

\usepackage{algorithm}
\usepackage[noend]{algpseudocode}
\usepackage{alphalph,etoolbox}

\usepackage{tikz}
\usepackage{tikz-cd}
\usepackage{pgfplots}
\usetikzlibrary{
  decorations.shapes,
  shapes.geometric,
  shadows,
  arrows,
  automata,
  positioning,
  calendar,
  mindmap,
  backgrounds,
  scopes,
  chains,
  er,
  3d,
  matrix,
  fit,
  intersections,
  pgfplots.fillbetween,
  pgfplots.statistics,
  pgfplots.colorbrewer
}
\usepgfplotslibrary{colorbrewer}
\pgfplotsset{compat=1.15}

\usepackage[colorlinks=true,citecolor=blue,urlcolor=blue]{hyperref}

\allowdisplaybreaks
\def\BibTeX{{\rm B\kern-.05em{\sc i\kern-.025em b}\kern-.08em
  T\kern-.1667em\lower.7ex\hbox{E}\kern-.125emX}}

\DeclareMathOperator*{\argmin}{argmin}
\DeclareMathOperator*{\argmax}{argmax}
\newcommand{\R}{{\rm I\!R}}
\newcommand{\N}{{\rm I\!N}}
\newcommand{\E}{{\rm I\!E}}
\newcommand{\T}{\mathbf{T}}
\newcommand{\A}{\mathcal{A}}
\newcommand{\AVaR}{\operatorname{AV@R}}
\newcommand{\EVaR}{\operatorname{EV@R}}

\newcommand{\trace}{\operatorname{tr}}
\newcommand{\conv}{\operatorname{conv}}
\newcommand{\sign}{\operatorname{sign}}
\newcommand{\prob}{\mathsf{P}}
\newcommand{\cov}{\mathrm{Cov}}
\newcommand{\var}{\mathrm{Var}}
\newcommand{\dd}{\mathrm{d}}
\newcommand{\lip}{\operatorname{lip}}
\newcommand{\diag}{\operatorname{diag}}
\newcommand{\eqd}{{}\overset{\rm d}{\sim}{}}

\definecolor{HeatA}{HTML}{FBFDFF}
\definecolor{HeatB}{HTML}{F1F7FC}
\definecolor{HeatC}{HTML}{E6F1F9}
\definecolor{HeatD}{HTML}{D8E9F5}
\definecolor{HeatE}{HTML}{C8DFF0}
\definecolor{HeatF}{HTML}{B5D3EA}
\newcommand{\heatA}[1]{\cellcolor{HeatA}\strut #1}
\newcommand{\heatB}[1]{\cellcolor{HeatB}\strut #1}
\newcommand{\heatC}[1]{\cellcolor{HeatC}\strut #1}
\newcommand{\heatD}[1]{\cellcolor{HeatD}\strut #1}
\newcommand{\heatE}[1]{\cellcolor{HeatE}\strut #1}
\newcommand{\heatF}[1]{\cellcolor{HeatF}\strut #1}
\newcommand{\heatlegendbox}[1]{%
  \begingroup\setlength{\fboxsep}{0pt}%
  \colorbox{#1}{\rule{1.15em}{0pt}\rule{0pt}{1.05ex}}%
  \endgroup}

\newtheorem{theorem}{Theorem}
\newtheorem{lemma}[theorem]{Lemma}
\newtheorem{proposition}[theorem]{Proposition}
\newtheorem{corollary}[theorem]{Corollary}
\newtheorem{assumption}{Assumption}
\newtheorem{remark}{Remark}

\usepackage{tcolorbox}
\tcbuselibrary{breakable}
\definecolor{propositionbg}{RGB}{255,250,230}

\newenvironment{proof}[1][Proof]{%
  \begin{IEEEproof}[#1]\ignorespaces
}{%
  \end{IEEEproof}
}

\usepackage{acronym}
\acrodef{aDP}[aDP]{approximate dynamic programming}
\acrodef{MLP}[MLP]{multilayer perceptron}
\acrodef{ICNN}[ICNN]{input convex neural network}
\acrodef{SDDP}[SDDP]{stochastic dual dynamic programming}
\acrodef{SOCP}[SOCP]{second order conic problem}
\acrodef{IQR}[IQR]{interquartile range}
\acrodef{pdf}[pdf]{probability density function}
\acrodef{cdf}[cdf]{cumulative distribution function}
\acrodef{LQ}[LQ]{linear-quadratic}

\newcommand{\paperfigure}[1]{\csname paperfigure@#1\endcsname}
\expandafter\newcommand\csname paperfigure@adversarial_cost_survival_h15_x2p25_1p5_a05\endcsname{%
  \includegraphics{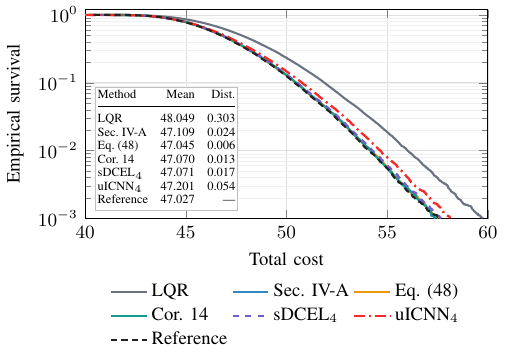}%
}
\expandafter\newcommand\csname paperfigure@adversarial_cost_survival_h15_x2p25_1p5_a09\endcsname{%
  \includegraphics{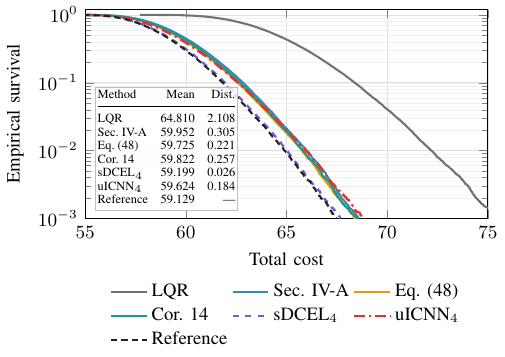}%
}
\expandafter\newcommand\csname paperfigure@adversarial_cost_survival_h15_common_lqr\endcsname{%
  \includegraphics{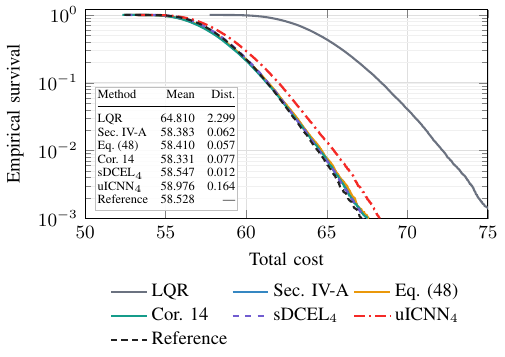}%
}
\expandafter\newcommand\csname paperfigure@adversarial_cost_survival_nx10_h10_energy3_a09\endcsname{%
  \includegraphics{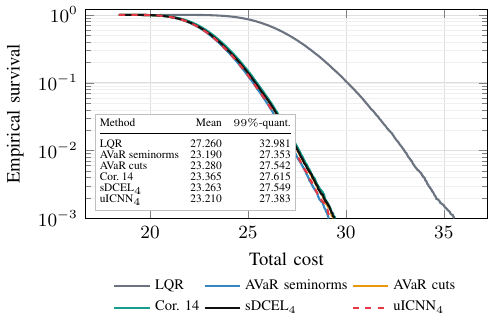}%
}
\expandafter\newcommand\csname paperfigure@bellman_rms_scaled_nx10_energy6p25\endcsname{%
  \includegraphics{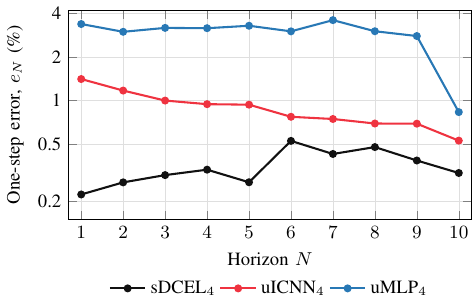}%
}

\begin{document}

\ifdefined\ArxivVersion
    \markboth{}{P. Sopasakis \& E. Vlahakis: Risk-Averse LQ Optimal Control}
\else
    \markboth{IEEE Transactions on Automatic Control}
    {P. Sopasakis \& E. Vlahakis: Risk-Averse LQ Optimal Control}
\fi

\title{Value function approximations for risk-averse linear-quadratic optimal control}

\author{Pantelis~Sopasakis and Eleftherios~E.~Vlahakis \IEEEmembership{Member, IEEE}%
\thanks{E. E. Vlahakis gratefully acknowledges financial support from KK-stiftelsen under Grant No. 202400088.}
\thanks{Pantelis Sopasakis is with the School of Electronics, Electrical Engineering and Computer Science, Queen's University Belfast, Northern Ireland, BT9 5BN, UK (e-mail: \texttt{p.sopasakis@qub.ac.uk}).}%
\thanks{Eleftherios E. Vlahakis is with the Department of Engineering Science, Division of Computer engineering and computer science, University West, 461 86 Trollhättan, Sweden (e-mail: \texttt{eleftherios.vlahakis@hv.se}).}%
}

\maketitle

%%%%%%%%%%%%%%%%%%%%%%%%%%%%%%%%%%%%%%%%%%%%%%%%%%%%%%%%%%%%%%%%%%%%%%%%%%%%%%%%
\begin{abstract}
Classical linear-quadratic control is analytically derivable because its value functions
remain quadratic and are characterized by a Riccati recursion. 
In multistage nested risk-averse optimal control, however, 
this quadratic structure is generally lost, even when
the dynamics and stage costs remain linear and quadratic. 
We study finite-horizon problems with additive symmetric
disturbances and law-invariant coherent risk measures. We show that each
optimal value function decomposes into the sum of a Riccati quadratic and a convex, 
even, globally Lipschitz residual. 
We exploit this structure to propose a seminorm-based value function 
approximation with global error bounds and a finite-horizon policy
suboptimality bound. We also develop tighter lower approximations using
supporting hyperplanes and propose specific computational recipes 
for the average value-at-risk and entropic value-at-risk. 
Numerical experiments demonstrate that the proposed approximation methods 
approximate the value functions accurately at a low computational cost.
Lastly, we leverage the value function decomposition to 
devise a structured approximate dynamic programming
method that retains the Riccati term and approximates only the residual.
Numerical experiments with Gaussian and exponential disturbances show that
the analytical approximations achieve high accuracy with low computational effort, 
while the proposed structure-preserving \ac{aDP} approach generally outperforms 
unstructured alternatives in value-function approximation.
\end{abstract}


\section{Introduction}
\label{sec:introduction}

\textbf{Motivation and problem statement.}
\Ac{LQ} optimal control is one of the fundamental benchmark problems 
in stochastic control, largely because of its analytical structure. 
For linear dynamics, quadratic costs, and additive zero-mean disturbances, 
the dynamic programming recursion preserves quadratic value functions, 
and the optimal feedback laws are characterized through a Riccati recursion. 
This has made \ac{LQ} control a classical setting and a building block for more 
general stochastic control problems. 

Risk-neutral stochastic optimal control 
typically minimizes the expected accumulated cost. However, the 
expectation may underrepresent low-probability realizations that incur large costs. 
A coherent risk measure provides a mechanism for assigning greater importance 
to such adverse outcomes while retaining useful convexity properties. 
Under some regularity assumptions, its dual representation allows the risk 
of a random cost to be interpreted as the worst-case expectation over a family of 
probability measures \cite{ShapiroDentcheva2021}. 
Conditional versions of such risk measures can be composed recursively 
to give rise to nested optimal control formulations, which are amenable
to a dynamic programming treatment \cite{RuszczynskiShapiro2006,Ruszczynski2010}.

This dual viewpoint gives risk-averse control a natural adversarial interpretation: 
at each stage, the controller selects a control action while an adversary selects 
an adverse probability distribution so as to maximize the 
expected cost. 
In this sense, each Bellman iteration can be viewed
as a sequential controller--adversary game.
Related game-theoretic interpretations have appeared 
in the distributionally robust control literature,
where uncertainty is represented through an opposing player acting on the disturbance law.
A recent example is the ``Devil'' formulation of \cite{SchobiLanzetti2026},
in which the adversary selects a disturbance distribution once, 
before the system evolution begins, and this distribution is then used throughout 
the evolution of the system.
The distinction with the present formulation is that, 
under the nested multistage risk criterion,
the adversarial distribution is selected separately 
at each stage based on the information available at that stage.
This way the adversary can adapt its choice to the current state 
and to the control action induced by the policy.
These stage-wise distributions arise endogenously from the dual representation 
of the coherent risk measure, while the nominal disturbance distribution remains fixed.

The loss of the quadratic structure of the value functions---and the 
corresponding loss of linearity of the optimal control laws---raises 
two key questions: what structure, if any, 
survives from classical \ac{LQ} control under nested coherent risk, 
and whether this structure be exploited to construct tractable controllers 
with quantitative approximation guarantees?

In this paper, we address these questions for discrete-time linear systems 
with quadratic costs and additive symmetric disturbances under law-invariant 
coherent risk measures. The key result is that the underlying \ac{LQ} structure 
is not entirely lost: 
the optimal value function at each stage admits the decomposition
\begin{equation}
    \label{eq:regularized-LQR}
    V_j^\star(x)
    {}={} 
    x^\intercal P_j x {}+{} r_j(x),
\end{equation}
where $P_j$ is generated by the Riccati recursion, 
while the residual $r_j$ is even, convex, and globally Lipschitz. 
The structure in \eqref{eq:regularized-LQR} provides a natural interpretation 
of risk-averse \ac{LQ} control as a regularized risk-neutral \ac{LQ} problem, 
where the classical Riccati term is retained and the residual acts as a 
convex regularizer.

This decomposition enables the development of approximation 
methods for the value functions. 
We firstly recursively construct an approximation of the residuals, $r_j$,
in terms of seminorms induced by the coherent risk measure 
which are convex, even, and globally Lipschitz. 
We show that the resulting value-function approximations 
satisfy global error bounds, and determine suboptimality bounds 
for the approximate feedback policies.
The approximate feedback laws are obtained from strongly convex 
one-step optimization problems; for Gaussian disturbances and 
average value-at-risk ($\AVaR$), these can be cast as
second-order-conic convex problems.

Next, we propose a tighter approximation method:
we sharpen this construction by exploiting the 
convexity of the risk-induced terms, $r_j$. 
The key idea is that through the dual representation of the 
coherent risk measure, the term $r_j$ can be written as the supremum 
of affine functions over a finite-dimensional image of the dual risk set. 
By selecting a finite subset of these affine functions 
we obtain a tighter approximation.
We give specific computational recipes for the cases of $\AVaR$ 
and the entropic value-at-risk ($\EVaR$).

Lastly, the structural characterization leads to a 
\textit{structured approximate dynamic programming} approach 
in which the known Riccati term is retained and only the convex, 
even, globally Lipschitz residual is learned from data.

\textbf{Related Work.} Risk-averse multistage optimization is 
closely connected to the theory of conditional risk mapping 
and dynamic risk measures. 
Conditional risk mappings and their recursive composition
are the building blocks of dynamically consistent risk-averse 
optimization and risk-averse dynamic programming 
\cite{RuszczynskiShapiro2006,Ruszczynski2010}.
Building up on this framework, Philpott \textit{et al.} 
\cite{PhilpottDeMatos2013}
develops a dynamic programming approach for multistage stochastic 
programs with coherent risk measures, using inner and outer 
approximations of the cost-to-go functions and incorporating
this into \ac{SDDP}.
Related work has studied time consistency, nested coherent risk measures,
and multistage stochastic optimization under nonlinear risk criteria \cite{FullnerRebennack2025,Pichler2022,Buerle2022}.

A closely related work is \cite{ChapmanLessard2022}
which studies a finite-horizon CVaR-LQ problem with distributional ambiguity and 
develops a dynamic-programming upper bound for its optimal value function. 
Instead this paper considers a nested risk formulation 
and a broader class of law-invariant coherent risk measures. 
More importantly, our starting point is an exact structural characterization 
of the risk-averse value function, from which certified approximations and 
feedback laws are subsequently derived.

A related body of work pertains to distributionally robust 
\ac{LQ} control \cite{shapiro_distributionally_2025}.
For example, \cite{Taskesen+2023} studies \ac{LQ} control when 
the disturbance distributions
are uncertain and are drawn from Wasserstein ambiguity sets,
and establishes conditions under which linear policies remain 
optimal.
Related formulations include distributionally robust 
infinite-horizon LQR \cite{hajar_distributionally_2024} and, 
more recently, Sinkhorn-based distributionally robust \ac{LQ} 
control \cite{cescon_global_2026}.
Game-theoretic interpretations of distributional 
uncertainty also arise in robust dynamic programming 
\cite{iyengar_robust_2005}.
In the recent work \cite{SchobiLanzetti2026}, 
the adversary responds to the chosen policy by fixing the 
disturbance laws over the horizon, while
our nested formulation instead allows the adversarial
disturbance law to depend on the current state and control
input at each stage.

Another relevant line of research concerns risk-averse multistage 
stochastic programming, where \ac{SDDP}
and related decomposition methods are used to approximate cost-to-go 
functions through successively generated affine cuts. These methods have 
been developed for multistage problems with linear costs and constraints 
under coherent risk measures 
\cite{PhilpottDeMatos2013,GuiguesRomisch2012,FreitasDaCosta2023}, 
including extensions that provide policy performance guarantees 
\cite{MerabetDaCostaLeclere2024} and applications to hydrothermal 
scheduling \cite{Diniz+2019}. More general formulations also extend 
these decomposition techniques to convex multistage optimization 
problems \cite{GuiguesMonteiro2021,FullnerRebennack2025}. 
SDDP formulations assume finitely supported disturbances,
with continuous distributions typically replaced by finite
sample-based approximations \cite[Sec.~11]{FullnerRebennack2025}.
Our supporting-hyperplane construction follows a similar 
geometric principle but exploits the specific structure of the 
risk-averse \ac{LQ} value function. Rather than approximating 
the entire cost-to-go function, we retain the Riccati
quadratic component exactly and use supporting hyperplanes to 
approximate an auxiliary risk-dependent function involved in 
constructing the nonquadratic residual. This leads to a 
structure-preserving approximation with global error guarantees.

Lastly, risk-averse extensions of \ac{aDP} \cite{Bertsekas2012_DPvol2_aDP} 
and  reinforcement learning have been developed for multistage 
stochastic decision problems, including simulation-based approaches 
\cite{JiangPowell2018}, 
risk-aware Markov decision processes \cite{YuHaskellXu2018}, 
and deep Q-learning under dynamic risk criteria \cite{YuShen2022}. 
Related work includes temporal-difference learning with 
linear value-function approximations under 
Markov coherent risk measures \cite{KoseRuszczynski2021}. 
These approaches employ numerical approximation methods 
without exploiting the specific analytical structure of 
risk-averse \ac{LQ} control. In contrast, our approach retains 
the Riccati quadratic component exactly and approximates 
only the residual. This analytical characterization additionally 
establishes convexity, evenness, and global Lipschitz continuity 
of the residual. Input-convex neural networks provide one mechanism 
for enforcing convexity \cite{AmosXiKotler2017}, while symmetrization 
and suitable norm constraints can enforce evenness and control 
Lipschitz continuity \cite{Gouk+2021}. 
This motivates the structured approximate dynamic programming 
approach investigated in this paper.

\textbf{Contributions.} The main contributions of the paper are summarized as follows:
(i) We characterize the exact value functions of nested risk-averse 
\ac{LQ} control as the sum of a Riccati-generated quadratic term and an even, 
convex, globally Lipschitz residual, and establish existence and 
uniqueness of the corresponding one-step optimal control action;
(ii) Building on this structure, we construct tractable approximations 
of the value functions and feedback laws, derive global uniform error bounds, 
and obtain an explicit finite-horizon suboptimality bound for 
the resulting approximate policy.
(iii) We develop tighter supporting-hyperplane approximations 
from the dual representation of coherent risk measures, 
including constructions for AV@R and EV@R, adaptive anchor selection, 
and analytical lower cuts for Gaussian disturbances. We also exploit the 
identified value-function structure in structured approximate dynamic programming.

\subsection{Notation}
Let $X$ and $Y$ be normed spaces; we say that a function 
$F:X \to Y$ is Lipschitz continuous
if there is $L \geq 0$ such that for all $x, x'\in X$
$\|F(x)-F(x')\| \leq L \|x-x'\|$.
The Lipschitz modulus of a function $F:\R^n\to\R^m$  is
$\lip F(\bar{x}) = \limsup_{x,x'\to \bar{x}, x\neq x'} \|F(x')-F(x)\| / \|x'-x\|$  \cite[Def. 9.1]{RockafellarWets2009}.
The sets of $n\times n$ symmetric positive (semi)definite matrices 
are denoted by ($\mathbb{S}_+^n$) $\mathbb{S}_{++}^n$.
For $Q\in\mathbb{S}_+^n$ we define the 
seminorm $\|x\|_Q^2 = x^\intercal Q x$.
The plus operator is defined as $(z)_+ = \max\{0, z\}$,
for $z\in\R$.
For a measurable set $E\in\mathcal{F}$ of a probability space 
$(\Omega, \mathcal{F}, \prob)$,
we define the random variable $1_E(\omega) = 1$, if $\omega\in \Omega$, 
and $0$ otherwise.
The \ac{pdf} and \ac{cdf} of the standard normal distribution,
$\mathcal{N}(0,1)$, are denoted by $\phi$ and $\Phi$,
respectively.
Two random variables, $Z_1$ and $Z_2$ being equal in 
distribution is denoted by $Z_1 \eqd Z_2$.
The support function of a nonempty set $C \subseteq \R^n$ is $\sigma_C(y) = \sup_{x\in C} y^\intercal x$.
The conjugate of a convex extended-real-valued function 
$f:\R^n\to\R\cup\{\infty\}$ is $f^*(y) {}={} \sup_x \{y^\intercal x - f(x)\}$.

\section{Problem statement}

\subsection{Risk measures}

A risk measure on the space $\mathcal{Z} \coloneqq \mathcal{L}^p(\Omega, \mathcal{F}, \prob; \R)$, $p\in [1, \infty)$ is a mapping 
$\rho:\mathcal{Z}\to\R\cup\{+\infty\}$. 
We regard random variables $Z\in \mathcal{Z}$ as costs, 
such that larger values of $Z$ are less desirable.
Roughly speaking, $\rho(Z)$ can quantify the size of the right tail 
of the distribution of $Z$.
A risk measure $\rho$ is called \emph{coherent} if it satisfies the following properties for all $Z, Z_1, Z_2 \in \mathcal Z$, $c \in \R$, and $\lambda > 0$ 
\cite[Sect. 6.3]{ShapiroDentcheva2021}
\begin{enumerate}[
    label=A\arabic*.,
    labelindent=1em,
    leftmargin=*
    ]
    % A1. Monotonicity
    \item
    \hypertarget{A1}{} (Monotonicity) 
    $\rho(Z_1)\leq \rho(Z_2)$, when $Z_1 \leq Z_2$ a.s.

    % A2. Translation equivariance
    \item
    \hypertarget{A2}{} (Translation equivariance) 
    $\rho(Z+c) = \rho(Z) + c$

    % A3. Positive homogeneity
    \item
    \hypertarget{A3}{} (Positive homogeneity) 
    $\rho(\lambda Z) = \lambda \rho(Z)$, and $\rho(0) = 0$

    % A4. Subadditivity
    \item 
    \hypertarget{A4}{} (Subadditivity)
    $\rho(Z_1 + Z_2) \leq \rho(Z_1)+\rho(Z_2)$
\end{enumerate}
Note that \hyperlink{A3}{A3} and \hyperlink{A4}{A4} \
imply that $\rho$ is convex, i.e.,
\begin{equation}
    \rho\bigl(\theta Z_1+(1-\theta)Z_2\bigr)
    \leq
    \theta\rho(Z_1)+(1-\theta)\rho(Z_2),
    \label{eq:risk_convexity}
\end{equation}
for $\theta\in[0,1]$. A coherent, proper, 
lower semicontinuous risk measure 
admits the dual representation 
\cite[Thm. 6.5]{ShapiroDentcheva2021}
\begin{equation}\label{eq:dual_rep_rho}
    \rho(Z) = \sup_{y \in \A}\E[y Z],
\end{equation}
where $\A$ 
is a nonemtpy convex subset of the dual space 
$\mathcal{Z}^* = \mathcal{L}^q(\Omega, \mathcal{F}, \prob)$
with $q\in (1, \infty]$ and $\nicefrac{1}{p} + \nicefrac{1}{q} = 1$ and for every $y\in \A$, 
$y \geq 0$ a.s. and $\E[y] = 1$.
We shall refer to $\A$ as the \textit{dual set}
of $\rho$.

For each $y\in\A$ we can write $\E[yZ]=\E^{\mu_y}[Z] = \int Z \dd \mu_y$,
where $\mu_y$ is the probability measure with 
$\dd \mu_y = y \dd \prob$ \cite{ShapiroDentcheva2021}; 
in particular, it is the measure
$\mu_y[B] = \E[y1_B(w)]$, over Borel sets $B$.
We can then express $\rho$ as a worst-case expectation 
over the probability measures of $\mathfrak{M} = \{\mu_y: y\in \A\}$, that is 
\begin{equation}
    \rho(Z) = \sup_{\mu\in \mathfrak{M}}\E^{\mu}[Z].
    \label{eq:dual_rep_rho_v2}
\end{equation}

A risk measure $\rho$ is called \emph{law-invariant} if
$\rho(Z_1) = \rho(Z_2)$ whenever 
$Z_1 \eqd Z_2$.
Hence, a law-invariant risk measure assigns the same risk to random costs having the
same distribution, irrespective of the underlying probability-space representation.

In this work, we consider law-invariant coherent 
risk measures, such as the popular average 
value-at-risk with parameter $\alpha\in [0, 1)$,
which is defined as 
\begin{equation}
    \AVaR_\alpha[Z] {}={} 
    \inf_{t\in\R}
    \bigl\{t + \tfrac{1}{1-\alpha}\E[(Z-t)_+]\bigr\}.
    \label{eq:avar_definition}
\end{equation}
If $Z$ is a continuous random variable, then $\AVaR_\alpha$ is the tail expectation
\begin{equation}
    \AVaR_\alpha[Z] 
    {}={} 
    \E\bigl[Z {}\mid{} Z \geq F_Z^{-1}(\alpha)\bigr],
    \label{eq:avar-tail-definition}
\end{equation}
where $F_Z^{-1}(\alpha)$ is the $\alpha$-quantile of $Z$. 
Note that for $\alpha=0$, $\AVaR_0[Z] = \E[Z]$.
The average value-at-risk of a normally distributed 
random variable can be computed based on the 
following proposition.

\begin{proposition}[$\AVaR$ of normal~\cite{RockafellarRoyset2010}]%
    \label{prop:avar_normal}
    Let $X\sim \mathcal{N}(0, 1)$
    and define $\kappa_\alpha \coloneqq \AVaR_\alpha[X]$.
    Then,
    \(
        \kappa_\alpha = \frac{\phi(\Phi^{-1}(\alpha))}{1-\alpha}
    \)
    for $\alpha \in (0, 1)$ and 
    $\kappa_0 = 0$.
\end{proposition}

Another risk measure is the entropic value-at-risk
with parameter $\alpha\in [0, 1)$, denoted $\EVaR_\alpha$, 
which is given by~\cite{Javid2012}
\begin{equation}
    \EVaR_\alpha[Z] 
    {}={}
    \inf_{t>0}
    \left\{
        t^{-1}\log \tfrac{\E[e^{tZ}]}{1-\alpha}
    \right\}. 
    \label{eq:evar-definition}
\end{equation}
The entropic value-at-risk of the standard normal distribution 
is given below.
\begin{proposition}[$\EVaR$ of normal \cite{Javid2012}]%
    \label{prop:evar_normal}
    Suppose $X\sim \mathcal{N}(0, 1)$.
    Then, for $\alpha \in (0, 1)$,
    $\EVaR_\alpha[X] = \sqrt{-2\log(1-\alpha)}$,
    and $\EVaR_{0}[X] = 0$.
\end{proposition}

\subsection{Risk-averse optimal control}

Consider the discrete-time linear system with 
an additive disturbance
\begin{equation}
    x_{k+1} = Ax_k + Bu_k + w_k,
    \label{eq:system-dynamics}
\end{equation}
with state $x_k\in\R^{n_x}$ and input $u_k\in\R^{n_u}$.
Hereafter we assume that the following assumption holds.
\begin{assumption}\label{assm:rho_disturbance}
    (i) The risk measure $\rho$ is coherent, proper,
    law-invariant, and lower semicontinuous,
    (ii) $w_k$ are iid random vectors
    with $w_k \overset{\rm d}{\sim} -w_k$, for all $k$,
    absolutely continuous with respect to Lebesgue measure, 
    and $\rho(\|w_k\|^2) < \infty$.
\end{assumption}

Consider the quadratic stage cost function 
\(
    \ell(x, u) = \|x\|_Q^2 + \|u\|_R^2,
\)
where $Q\in\mathbb{S}_{+}^{n_x}$ and 
$R\in\mathbb{S}_{++}^{n_u}$, and the terminal cost 
\(
    V_f(x) = \|x\|_{P_f}^2,
\)
with $P_f\in\mathbb{S}_{+}^{n_x}$. The objective is 
to determine a finite sequence of Borel control laws, 
$\bm{\pi}_N = (\pi_N, \ldots, \pi_1)$,
via the multistage risk-averse optimal control problem 
\cite{Ruszczynski2010}
\begin{subequations}
    \label{eq:risk-averse-problem-original}
    \begin{align}
        \operatorname*{Minimize}_{\bm{\pi}_N \in \Pi_N}
        \;
        &\ell(x_0,\pi_N(x_0))
        +\rho_{w_0}\Bigl[
        \ell(x_1,\pi_{N-1}(x_1))
        \notag\\
        &
        +\cdots+
        \rho_{w_{N-2}}\bigl[
            \ell(x_{N-1},\pi_1(x_{N-1}))
        \notag\\
        & \quad +\rho_{w_{N-1}}[V_f(x_N)]
        \bigr]\cdots
        \Bigr],
        \\
        \text{subject to:}\;& x_{k+1} = Ax_k + B\pi_{N-k}(x_k) + w_k,
        \notag\\
        &\quad k=0,\ldots, N-1,
    \end{align}
\end{subequations}
with horizon $N$,
for a fixed initial state $x_0$, 
and where $\rho_{w_k}$ is a risk measure with respect 
to the random variable $w_k$ with $w_0, \ldots, w_{k-1}$
held fixed.
The minimization in \eqref{eq:risk-averse-problem-original}
is over sequences of Borel feedback laws,
$\pi_k:\R^{n_x}\to\R^{n_u}$, 
$k=1,\ldots, N$,
such that the nested objective is finite. 
The corresponding dynamic programming recursion is
\begin{subequations}%
    \begin{align}%
        V_{j+1}^\star(x)
        {}={}& 
        (\T V_j^\star)(x),
        \label{eq:dp-value-iteration}
        \\
        \kappa_{j+1}^\star(x) 
        {}\in{}& 
        (\mathbf{S}V_j^\star)(x),
    \end{align}
\end{subequations}
with $V_0^\star(x) {}={} V_f(x)$, where 
\begin{subequations}
    \begin{align}
        (\T V)(x) 
        \coloneqq{}&
        \inf_{u}\bigl\{\ell(x, u) + \rho [V(Ax + Bu + w)]\bigr\},
        \\
        (\mathbf{S}V)(x) 
        \coloneqq{}&
        \argmin_{u}\bigl\{\ell(x, u) {+} \rho [V(Ax {+} Bu {+} w)]\bigr\}.
    \end{align}
\end{subequations}

In light of \eqref{eq:dual_rep_rho_v2}, each dynamic 
programming iteration admits a Stackelberg game interpretation: 
for a given state $x$, the controller first chooses $u$, 
and the cost-to-go is then evaluated by taking the 
supremum of the expected next-stage cost-to-go, 
$\E^\mu[V_j(Ax+Bu+w)]$, over $\mu\in\mathfrak{M}$. 
This adversarial viewpoint is reminiscent of ``the Devil'' 
in \cite{SchobiLanzetti2026}, although in that paper 
the adversary selects disturbance laws for the 
entire horizon after observing the control policy, 
whereas here the adversarial law may depend on the 
stage, state, and control. The difference between the 
two formulations is discussed in \cite[Sec. 4]{Shapiro2012}.

For a function $V:\R^{n_x}\to\R$ and a controller $\pi:\R^{n_x}\to\R^{n_u}$ 
we define
\begin{equation}
    (\T_\pi V)(x) = \ell(x, \pi(x)) + \rho [V(Ax + B\pi(x) + w)],
\end{equation}
provided $V$ and $\pi$ satisfy Borel-measurability and 
integrability assumptions for the risk to be well-defined 
and finite. 
For a policy $\bm{\pi}_N=(\pi_{N}, \ldots, \pi_{1})$ we denote 
$\T_{\bm{\pi}_N} = \T_{\pi_{N}}\ldots \T_{\pi_{1}}$ and define 
$J_{\bm{\pi}_N}(x) = (\T_{\bm{\pi}_N} V_0^\star)(x)$.
From the dynamic programming recursion, the problem in 
\eqref{eq:risk-averse-problem-original} can be written as 
\begin{equation}
    \operatorname*{Minimize}_{\bm{\pi}_N \in \Pi_N} 
    J_{\bm{\pi}_N}(x).
\end{equation}

% From Proposition \ref{prop:avar_normal},
% $\AVaR_\alpha[\mathcal{N}(0, \sigma^2)]=\sigma \kappa_\alpha$,
% and for $w\sim \mathcal{N}(0, \Sigma_w)$ and $c\in\R^{n_x}$,  $\AVaR_\alpha[c^\intercal w] 
%     {}={}         
%     \|\Sigma_w^{1/2}c\|_2 \kappa_\alpha$.

\section{Structure of the value functions}
Unlike in the deterministic or expectation-based linear-quadratic 
case, functions $V_j^\star$ are not quadratic.  
However, they possess the properties stated in the following 
theorem: $V_j^\star$ can be decomposed into a sum of the 
LQR cost function $\|x\|_{P_j}^2$ and a convex, even, and 
globally Lipschitz residual $r_j(x)$.
\begin{theorem}[Structure of $V_j^\star$]
    \label{thm:Vjstar_structure}
    For all $j$, functions $V_j^\star$ are convex and have the form 
    \begin{equation}
        V_j^\star(x) = \|x\|_{P_j}^2 + r_j(x),
        \label{eq:Vjstar_structure}
    \end{equation}
    where the matrices $P_j$ are generated by a 
    Riccati recursion and $r_j$ are even, convex,
    and globally Lipschitz functions.
    Moreover, for every $x$, the infimum defining
    $(\T V_j^\star)(x)$ is attained and
    $(\mathbf{S}V_j^\star)(x)$ is a singleton.
\end{theorem}

Theorem \ref{thm:Vjstar_structure} shows that the 
nested risk-averse value function retains the 
Riccati quadratic of standard LQR, with an 
additional risk-induced term. 
This way the problem 
can be viewed as a \textit{regularized LQR problem}.
Moreover, the Riccati-based decomposition in 
\eqref{eq:Vjstar_structure} motivates the development 
of methodologies for the approximation of the residual,
$r_j$.

To prove this theorem we first need the following  
lemma 

\begin{lemma}[Properties of $\psi$]\label{lem:psi_properties}
    Let $P\in\mathbb{S}_{+}^{n_x}$, let
    $r:\R^{n_x}\to\R$ be convex, even, and globally
    $L_r$-Lipschitz.
    Define $\psi:\R^{n_x}\to\R$ by
    \begin{equation}
        \psi(m)=\rho
        \bigl[2m^\intercal Pw+\|w\|_P^2+r(m+w)\bigr].
    \end{equation}
    Then, $\psi$ is finite-valued, convex, even, and globally
    Lipschitz.
\end{lemma}

The proof of Lemma \ref{lem:psi_properties} can be found in 
the \hyperlink{link:appendix}{Appendix}.
We can now prove Theorem \ref{thm:Vjstar_structure}.

\begin{proof}[Proof of Theorem~\ref{thm:Vjstar_structure}]
The result holds at $j=0$ with $P_0=P_f$ and $r_0=0$.
Suppose it holds at $j$ and set $m=Ax+Bu$. Translational
equivariance of $\rho$ gives
\begin{align}
    (\T V_j^\star)(x)
    =\inf_u\{&\|x\|_Q^2+\|u\|_R^2+\|m\|_{P_j}^2
        +\psi_j(m)\},
\end{align}
where
\begin{equation}
    \label{eq:psij}
    \psi_j(m)=\rho
    [2m^\intercal P_jw+\|w\|_{P_j}^2+r_j(m+w)].
\end{equation}
By Lemma~\ref{lem:psi_properties}, $\psi_j$ is finite-valued,
convex, even, and globally Lipschitz.

Let us define $\bar{R}_j=R+B^\intercal P_jB$,
the gain matrices $K_{j+1}=-\bar{R}_j^{-1}B^\intercal P_jA$,
the closed-loop matrices $\bar{A}_{j+1}=A+BK_{j+1}$, and 
the recursion
\begin{equation}
    P_{j+1}=Q+A^\intercal P_jA-A^\intercal P_jB\bar{R}_j^{-1}B^\intercal P_jA,
\end{equation}
for $j=0, \ldots, N-1$.
Expanding the quadratic terms and using these definitions gives
$\|x\|_Q^2+\|u\|_R^2+\|Ax+Bu\|_{P_j}^2 = \|x\|_{P_{j+1}}^2 +\|u-K_{j+1}x\|_{\bar{R}_j}^2$, therefore
\begin{align}
    (\T V_j^\star)(x)
    =\inf_u\{&\|x\|_{P_{j+1}}^2
      +\|u-K_{j+1}x\|_{\bar{R}_j}^2 \notag\\
      &+\psi_j(Ax+Bu)\}.
      \label{eq:TVjstar_psi}
\end{align}
Let us introduce the auxiliary input $v=u-K_{j+1}x$. Since $Ax+Bu=\bar{A}_{j+1}x+Bv$, we have
\begin{align}
    (\T V_j^\star)(x)
    {}={}&
    \|x\|_{P_{j+1}}^2+r_{j+1}(x), 
\intertext{with}
    r_{j+1}(x)
    {}\coloneqq{}&
    \inf_v
      \{\|v\|_{\bar{R}_j}^2+\psi_j(\bar{A}_{j+1}x+Bv)\}.
      \label{eq:r-j-plus-one}
\end{align}
Moreover, note that $P_{j+1}\succeq0$.

For fixed $x$, denote the objective defining $r_{j+1}(x)$ by
$g_x(v)$ and let $L_{\psi_j}$ be a Lipschitz constant of $\psi_j$.
Then
\begin{equation*}
    g_x(v)\geq\lambda_{\min}(\bar{R}_j)\|v\|^2
      +\psi_j(\bar{A}_{j+1}x)-L_{\psi_j}\|B\|_2\|v\|,
\end{equation*}
so $g_x$ is coercive; it is also continuous and strongly convex,
thus it has a unique minimizer; the infimum above is a
minimum, and $u{=}K_{j+1}x{+}v$ is the unique element of
$(\mathbf S V_j^\star)(x)$.

The function $r_{j+1}$ is convex as the inf-projection
of a convex function \cite[Prop. 2.22]{RockafellarWets2009}. 
By substituting
$z=-v$ and using the evenness of $\psi_j$, one also obtains
$r_{j+1}(-x)=r_{j+1}(x)$. 
Lastly, for each $v\in\R^{n_u}$, define
$f_v(x)=g_x(v)$, so that
$r_{j+1}(x)=\inf_{v} f_v(x)$.
Since the mapping $x\mapsto\bar A_{j+1}x+Bv$ has Lipschitz
modulus $\|\bar A_{j+1}\|_2$, independent of $v$,
\cite[Ex.~9.8(c)]{RockafellarWets2009} gives
$\lip f_v(x) \leq L_{\psi_j}\|\bar A_{j+1}\|_2 \eqqcolon C_j$.
Applying \cite[Prop.~9.10]{RockafellarWets2009} to the
family $\{f_v\}_{v\in\R^{n_u}}$, with the 
upper-semicontinuous majorant $k(x) = C_j$, yields
$\lip r_{j+1}(x)\leq C_j$.
Since $r_{j+1}$ is defined over $\R^{n_x}$,
it is globally $C_j$-Lipschitz as it follows from 
\cite[Thm.~9.2]{RockafellarWets2009}.
\end{proof}

% \lef{
% \begin{remark}
%     Theorem~\ref{thm:Vjstar_structure} admits an interpretation
% of the risk-averse problem as a regularized LQR problem.
% Indeed, the optimal control law consists of the classical
% LQR feedback and a correction induced by the convex,
% even, and globally Lipschitz function $\psi_j$, which
% acts as a risk-dependent regularizer on the nominal
% successor state. In the risk-neutral case, this
% regularizer is constant, recovering the classical
% LQR feedback law.
% \end{remark}
% }

A direct consequence of Theorem \ref{thm:Vjstar_structure}
is that 
\begin{equation}
    V_j^\star(x)-V_j^\star(0)
    {}={}
    \|x\|_{P_j}^{2}
    {}+{}
    \bigl(r_j(x)-r_j(0)\bigr),
\end{equation}
where the residual satisfies 
\begin{equation}
    0\leq r_j(x)-r_j(0)\leq L_{r_j}\|x\|,
\end{equation}
that is, the centered residual $r_j(x)-r_j(0)$ grows 
at most linearly with $\|x\|$. 
This motivates approximating the residual while retaining 
the quadratic term in \eqref{eq:Vjstar_structure}.

Next we will see how the decomposition of 
Theorem \ref{thm:Vjstar_structure} can be used to 
produce bounded-error approximations of $V_j^\star$.

\section{Bounded value function approximations}%
\label{section:value-funct-approx}
Although $V_j^\star$ has the form given in \eqref{eq:Vjstar_structure},
functions $r_j$ cannot be computed explicitly. In this section we will 
approximate $V_j^\star$ by functions of the form 
$\widehat{V}_j(x) = \|x\|_{P_j}^2 + s_j(x)$,
and will establish approximation bounds. Our goal is 
to determine value function approximations, $\widehat{V}_j$,
and corresponding controllers, $\widehat{\pi}_{j+1}$, 
satisfying the following desiderata:
\begin{enumerate}[
  label=D\arabic*.,
  labelindent=1em,
  leftmargin=*
]
    \item \hypertarget{D1}{} There is $\epsilon_j > 0$ such that $\|\widehat{V}_j - V_j^\star\|_\infty \leq \epsilon_j$,
    \item \hypertarget{D2}{} There is $\delta_j > 0$ such that $0 \leq (\T_{\widehat{\pi}_{j+1}} \widehat{V}_j)(x) - (\T \widehat{V}_j)(x) \leq  \delta_j$, for all $x\in\R^{n_x}$.
\end{enumerate}

Under \hyperlink{D1}{D1} and \hyperlink{D2}{D2}, 
% $\leq J_{\widehat{\bm{\pi}}_N}(x)$
$J_{\widehat{\bm{\pi}}_N}(x)$ is a bounded-error approximation 
of $V_N^\star(x)$ as shown in the following proposition.

\begin{proposition}
    \label{prop:desiderata-consequences}
    If \hyperlink{D1}{D1} and \hyperlink{D2}{D2} are satisfied, then the 
    approximate policy 
    $\widehat{\bm{\pi}}_N = (\widehat{\pi}_N, \ldots, \widehat{\pi}_1)$
    satisfies 
    \begin{equation}
        0 \leq J_{\widehat{\bm{\pi}}_N}(x) - V_N^\star(x) \leq \sum_{j=0}^{N-1}(\delta_j + 2\epsilon_j),
        \label{eq:JN-VN-bound}
    \end{equation}
    for all $x\in\R^{n_x}$.
\end{proposition}

\begin{proof}
    If $\|V-V'\|_\infty\leq\epsilon$, monotonicity and translation
    equivariance of $\rho$ give
    $|\T_\pi V(x)-\T_\pi V'(x)|\leq\epsilon$ for every controller
    $\pi$ and every $x$. Taking infima over control actions 
    gives $|\T V(x)-\T V'(x)|\leq\epsilon$. 
    \hyperlink{D1}{D1} and \hyperlink{D2}{D2} imply
    \begin{align}
        0
        &\leq (\T_{\widehat\pi_{j+1}}V_j^\star)(x)
              -(\T V_j^\star)(x)
        \notag\\
        &\leq (\T_{\widehat\pi_{j+1}}\widehat V_j)(x)
              -(\T\widehat V_j)(x)+2\epsilon_j
        \leq \delta_j+2\epsilon_j.
        \label{eq:desiderata-one-step-bound}
    \end{align}
    The first inequality holds because $\T$ takes the 
    infimum over control actions.

    Let $J_0=V_0^\star$ and
    $J_{j+1}=\T_{\widehat\pi_{j+1}}J_j$, and set
    $S_j=\sum_{i=0}^{j-1}(\delta_i+2\epsilon_i)$, with $S_0=0$.
    Suppose that $0\leq J_j-V_j^\star\leq S_j$ pointwise.
    Monotonicity and translation equivariance give
    \begin{equation*}
        (\T_{\widehat\pi_{j+1}}V_j^\star)(x)
        \leq J_{j+1}(x)
        \leq (\T_{\widehat\pi_{j+1}}V_j^\star)(x)+S_j.
    \end{equation*}
    From \eqref{eq:desiderata-one-step-bound},
    $0\leq J_{j+1}(x)-V_{j+1}^\star(x)
    \leq S_j+\delta_j+2\epsilon_j=S_{j+1}$.
    Since $J_0=V_0^\star$, inductively \eqref{eq:JN-VN-bound} 
    holds because $J_N=J_{\widehat{\bm{\pi}}_N}$.
\end{proof}

Next, we consider approximations of the form 
\begin{equation}
    \widehat{V}_j(x) = \|x\|_{P_j}^2 + s_j(x),
    \label{eq:Vj_structured_approximation}  
\end{equation}
where $s_j$ are convex, even, and Lipschitz functions.
Our objective is to establish approximations of $\psi_j$ 
in \eqref{eq:psij} with uniformly bounded error. 
To this end, we turn our attention to the quantity
$\xi_{j}(m) = \|w\|_{P_j}^2 + 2m^\intercal P_j w$
and its risk, that is,
\begin{equation}
    \eta_j(m) {}\coloneqq{} \rho\left(\xi_j(m)\right).
    \label{eq:eta_j:definition}
\end{equation}
We will show that by appropriately approximating $\eta_j$
we can construct functions $s_j$ which guaratee the 
satisfaction of \hyperlink{D1}{D1} and \hyperlink{D2}{D2}.
In Section \ref{sec:sum-seminorms} we propose a simple 
approximation of $s_j$ by a sum of seminorms, while in 
Section \ref{sec:finer-approximation} we propose a tighter 
approximation using supporting hyperplanes.

Lastly, \hyperlink{D1}{D1} constitutes a global error bound 
on the approximation of $V_j^\star$, but 
\eqref{eq:Vj_structured_approximation} yields a bound 
near zero as well. For the approximate and exact value 
functions let us define 
$\Delta \widehat{V}_j(x) = \widehat{V}_j(x) - \widehat{V}_j(0)$, 
$\Delta V_j^{\star}(x) = V_j^{\star}(x) - V_j^{\star}(0)$ and 
the \textit{centered error}
\begin{align}
    e_j(x) {}={}& \Delta \widehat{V}_j(x) - \Delta V_j^{\star}(x)
    \\
    {}={}& \Delta s_j(x) - \Delta r_j(x),
\end{align}
where $\Delta s_j(x) = s_j(x) - s_j(0)$ and 
$\Delta r_j(x) = r_j(x) - r_j(0)$.
Since both $s_j$ and $r_j$ are convex, even, and Lipschitz,
it follows that 
\begin{equation}
    |e_j(x)| \leq \max\{L_{s_j}, L_{r_j}\}\|x\|,
    \label{eq:centered-error-Ox-bound}
\end{equation}
where $L_{s_j}$ is a Lipschitz constant of $s_j$.
Later we will show that this $O(\|x\|)$-error bound 
can be improved to an $O(\|x\|^2)$, locally, 
under certain conditions
(see Proposition \ref{prop:quadratic-centered-error}).

\subsection{Sum-of-seminorms approximation}\label{sec:sum-seminorms}
In this section we take $s_j(x)$ 
to be a sum of seminorms, which are convex, even, and 
Lipschitz.

\begin{proposition}[Uniform bounds for $\psi_j$]%
    \label{prop:psi_approx}
    For each $j=0,\ldots,N-1$, let
    $\tilde{\eta}_j(m) = 2\rho\left(m^\intercal P_j w\right)$.
    Then,
    \begin{equation}
        0
        \leq{}
        \eta_j(m) - \tilde{\eta}_j(m) 
        {}\leq{}
        \rho\bigl(\|w\|_{P_j}^2\bigr).
        \label{eq:avar_xi_bounds}
    \end{equation}
    Moreover, there is $C_j \geq 0$ such that 
    $|\psi_j(m) - r_j(m) - \tilde{\eta}_j(m)| \leq C_j$.
\end{proposition}

\begin{proof}    
    We have $\|w\|_{P_j}^2 \geq 0$, a.s., so 
    $\xi_j(m) \geq 2m^\intercal P_j w$,
    and by the monotonicity of $\rho$,
    the lower bound in \eqref{eq:avar_xi_bounds}
    is proven.
    By subadditivity and positive homogeneity, we obtain the 
    upper bound in \eqref{eq:avar_xi_bounds}.

    Let $L_{r_j}$ be a Lipschitz constant of $r_j$ 
    (Theorem \ref{thm:Vjstar_structure}), so 
    $|r_j(m+w)-r_j(m)|\leq L_{r_j}\|w\|$.
    Using \eqref{eq:rho_useful_ineq}, the translational equivariance of 
    $\rho$, the triangle inequality, and the Lipschitzianity of $r_j$,
    \begin{align}
        &|\psi_j(m) - r_j(m) - \tilde{\eta}_j(m)|
        \notag\\
        &{}={}
        |\rho\left(\xi_j(m)+r_j(m+w)\right)- r_j(m)-\tilde{\eta}_j(m)|
        \notag\\
        &{}={}
        |\rho\left(\xi_j(m)+r_j(m+w)- r_j(m)\right)-\tilde{\eta}_j(m)|
        \notag\\
        &{}\leq{}
        |\rho\left(\xi_j(m)+r_j(m+w)- r_j(m)\right) 
        % \notag\\
        - \rho\left(\xi_j(m)\right)| \notag\\
        &{}\quad{}
        + |\rho\left(\xi_j(m)\right) -\tilde{\eta}_j(m)|
        \notag\\
        &{}\leq{}
        \rho\left(|r_j(m+w) - r_j(m)|\right) + \rho(\|w\|_{P_j}^2) 
        \notag\\
        &{}\leq{} 
        L_{r_j}\rho\left(\|w\|\right) + \rho(\|w\|_{P_j}^2) \eqqcolon C_j.
        \label{eq:psi_j_bound}
    \end{align}
    By Assumption \ref{assm:rho_disturbance}, 
    this bound is finite
    and independent of $m$, proving the proposition.
\end{proof}

Proposition \ref{prop:psi_approx} allows us to establish a uniform 
approximation of $r_{j}$.
As a remark, if $w_k \sim \mathcal{N}(0, \Sigma_w)$, and $\rho = \AVaR_\alpha$ with $\alpha\in[0,1)$,
then $2m^\intercal P_j w \sim \mathcal{N}(0, 4m^\intercal P_j \Sigma_w P_j m)$, 
and from Proposition \ref{prop:avar_normal},
$\tilde{\eta}_j(m)=\AVaR_\alpha[2m^\intercal P_j w]=2\kappa_\alpha\|\Sigma_w^{1/2}P_jm\|$.
Note also that under Assumption \ref{assm:rho_disturbance}, $\tilde{\eta}_j$
is a seminorm for any law-invariant coherent 
risk measure $\rho$.

\begin{proposition}[Uniform bounds for $r_j$]%
    \label{prop:r_approx} For each $j=0,\ldots,N-1$, there is $D_j \geq 0$
    such that 
    $0\leq \psi_j(\bar{A}_{j+1}x)-r_{j+1}(x) \leq D_j$.
\end{proposition}

\begin{proof}
    Choosing $v=0$ in \eqref{eq:r-j-plus-one} gives
    \begin{equation}
        r_{j+1}(x)\leq \psi_j(\bar{A}_{j+1}x),
    \end{equation}
    for $x\in\R^{n_x}$. 
    Let $L_{\psi_j}$ be a Lipschitz constant of $\psi_j$
    (Lemma \ref{lem:psi_properties}). Then,
    $|\psi_j(\bar{A}_{j+1}x+Bv) - \psi_j(\bar{A}_{j+1}x)| \leq L_{\psi_j}\|B\|\;\|v\|$, so
    \begin{align}
        &\|v\|_{\bar{R}_j}^2+\psi_j(\bar{A}_{j+1}x+Bv) 
        \notag\\
        &\quad\geq \psi_j(\bar{A}_{j+1}x)
          +\lambda_{\min}(\bar{R}_j)\|v\|^2
          -L_{\psi_j}\|B\|_2\|v\| 
        \notag\\
        &\quad\geq \psi_j(\bar{A}_{j+1}x)
          -\frac{L_{\psi_j}^2\|B\|_2^2}
          {4\lambda_{\min}(\bar{R}_j)},
        \label{eq:prop:r_approx:lower_bound}
    \end{align}
    where the last inequality follows by minimizing the scalar quadratic
    $\lambda_{\min}(\bar{R}_j)t^2-L_{\psi_j}\|B\|_2t$ over $t\geq0$.
    Taking the infimum over $v$ on both sides of \eqref{eq:prop:r_approx:lower_bound},
    \begin{equation*}
        0\leq \psi_j(\bar{A}_{j+1}x)-r_{j+1}(x)
        \leq \frac{L_{\psi_j}^2\|B\|_2^2}
        {4\lambda_{\min}(\bar{R}_j)}
        {}\eqqcolon{} D_j.
    \end{equation*}
    This completes the proof.
\end{proof}

We now define the sequence of functions 
$s_0(x)=0$ and 
\begin{equation}
    s_{j+1}(x)=s_j(\bar{A}_{j+1} x) +\tilde{\eta}_j(\bar{A}_{j+1} x).
    \label{eq:dfn_s}
\end{equation}
We can see that 
\begin{equation}
    s_j(x)
    {}={}
    \sum_{i=0}^{j-1}
    \tilde{\eta}_i\bigl(\bar{A}_{i+1}\cdots \bar{A}_{j-1}\bar{A}_{j}x\bigr).
\end{equation}

Replacing $\psi_j(m)$ by
$s_j(m)+\tilde{\eta}_j(m)$ in \eqref{eq:TVjstar_psi}
we can define the approximate controller
\begin{align}
    \widehat{\kappa}_{j+1}(x)
    {}\coloneqq{}&
    \hspace{-0.5em}
    \argmin_{
        \substack{
            u\in\R^{n_u} \\ m = Ax+Bu
        }
        }
    \hspace{-0.3em}
    \|u\|_R^2
    + \|m\|_{P_j}^2
    + s_j(m) + \tilde{\eta}_j(m)
    \notag\\
    {}={}&
    K_{j+1}x {}+{}
    \hspace{-2.5em}
    \argmin_{
        \substack{
            v\in\R^{n_u} 
            \\ 
            z_{i} = \bar{A}_{i+1}z_{i+1}, 
            i=j-1,\ldots, 0
            \\
            z_j = \bar{A}_{j+1}x + Bv
        }
        }
    \hspace{-1.5em}
    \|v\|_{\bar{R}_j}^2
    + \sum_{i=0}^{j} \tilde{\eta}_i(z_i).
    \label{eq:kappa_hat_jp1}
\end{align}
The minimizer is unique because the objective is strongly convex in $u$.
Note that the $\widehat{\kappa}_{j+1}$ are determined by 
solving a strongly convex optimization problem. In the 
Gaussian case with the $\AVaR$ risk, 
these problems are \acp{SOCP}.
It remains to confirm that \hyperlink{D1}{D1} and \hyperlink{D2}{D2} 
are satisfied, so that by Proposition \ref{prop:desiderata-consequences}
we have an approximation of $\widehat{V}_N$.

\begin{proposition}
    \label{prop:sum-seminorms-desiderata}
    The sum-of-seminorms approximation satisfies 
    desiderata \hyperlink{D1}{D1} and \hyperlink{D2}{D2}.
\end{proposition}

\begin{proof}
    % Part 1: prove D1
    Set $z=\bar A_{j+1}x$. From Proposition \ref{prop:psi_approx} and 
    \ref{prop:r_approx} and by the triangle 
    inequality
    \begin{align*}
       \left|r_{j+1}(x)-s_{j+1}(x)\right|
       {}\leq{}&
       \left|r_{j+1}(x)-\psi_j(z)\right|
       \notag\\
       &{}+{}
       \left|\psi_j(z) - r_j(z) - \tilde\eta_j(z)\right|
       \notag\\
       &{}+{}
       \left|r_{j}(z) - s_j(z)\right|
       \notag\\
       {}\leq{}& 
       C_j + D_j + \left|r_{j}(z) - s_j(z)\right|.
    \end{align*}
    Since $r_0=s_0=0$, by induction
    \begin{equation}
        \|\widehat V_j-V_j^\star\|_\infty
        =\|s_j-r_j\|_\infty
        \leq\sum_{i=0}^{j-1}
             (C_i + D_i),
    \end{equation}
    so \hyperlink{D1}{D1} holds.

    % Part 2: prove D2
    To prove \hyperlink{D2}{D2}, let $L_{s_j}$ be a global
    Lipschitz constant of $s_j$. 
    For brevity let us define 
    $\bar{w} = \rho(\|w\|)$, and 
    $\tilde{w}_j = \rho(\|w\|_{P_j}^2)$.
    % Such a constant exists because
    % \eqref{eq:rho_useful_ineq} gives
    % $|\tilde\eta_i(m)-\tilde\eta_i(m')|
    % \leq 2\|P_i\|_2\AVaR_\alpha[\|w\|]\|m-m'\|$,
    % and \eqref{eq:dfn_s} involves only finite sums and
    % linear compositions.
    For $m\in\R^{n_x}$, define
    \begin{equation*}
        e_j(m)=
        \rho[\xi_j(m,w) + s_j(m{+}w)]
        -s_j(m)-\tilde\eta_j(m).
    \end{equation*}
    We then have
    % Translation equivariance and \eqref{eq:rho_useful_ineq} give
    \begin{align}
        &\bigl|\rho[\xi_j(m,w)+s_j(m+w)]
              -s_j(m)-\eta_j(m)\bigr|
        \notag\\
        &\quad
        {}={}
        \bigl|\rho\left[\xi_j(m,w)+s_j(m+w)-s_j(m)\right]
              -\rho(\xi_j(m))\bigr|
        \notag\\
        &\quad
        {}\overset{\eqref{eq:rho_useful_ineq}}{\leq}{}
        \rho\left[|s_j(m+w)-s_j(m)|\right]
        \leq L_{s_j}\bar{w}.
        \label{eq:ej_with_eta_bound}
    \end{align}
    We combine the bounds of \eqref{eq:avar_xi_bounds} with 
    \eqref{eq:ej_with_eta_bound} to obtain
    $-L_{s_j}\bar{w}
        \leq e_j(m)
        \leq \tilde{w}_j
             +L_{s_j}\bar{w}$.
    For $m=Ax+Bu$, let
    $h_x(u)=\ell(x,u)+\|m\|_{P_j}^2+s_j(m)+\tilde\eta_j(m)$;
    the controller \eqref{eq:kappa_hat_jp1} minimizes $h_x$.
    The objective defining
    $(\T\widehat V_j)(x)$ is $F_x(u)=h_x(u)+e_j(Ax+Bu)$. 
    Since $h_x(\widehat\kappa_{j+1}(x))\leq h_x(u)$ 
    for all $u$,
    \begin{align}
        F_x(\widehat\kappa_{j+1}(x))-F_x(u)
        &\leq \tilde{w}_j
        +2L_{s_j}\bar{w} \eqqcolon \delta_j.
    \end{align}
    Taking the infimum over $u$ gives the upper bound in
    \hyperlink{D2}{D2}.
\end{proof}

\subsection{A finer approximation via supporting hyperplanes}\label{sec:finer-approximation}
In Proposition \ref{prop:psi_approx} we established that
$\tilde{\eta}_j$ is a bounded-error approximation of $\eta_j$,
where $\tilde{\eta}_j$ is a seminorm.
A finer approximation can be obtained by exploiting 
the fact that $\eta_j$ is finite for all $m$ and 
convex, thus at every $m$ a subgradient exists~\cite[Thm. 23.4]{Rockafellar1970}.
Given a point $\hat{m}$ and $\hat{a}\in \partial \eta_j(\hat{m})$, 
let $\hat{b} = \eta_j(\hat{m}) - \hat{a}^\intercal \hat{m}$;
then $\hat{a}^\intercal m + \hat{b}$ defines a supporting 
hyperplane of $\eta_j$ at $\hat{m}$.
By the evenness of $\eta_j$, both 
$\pm\hat{a}^\intercal m + \hat{b}$ are supporting hyperplanes
(at $\pm \hat{m}$).
Using a collection of such anchor points, 
$(m_{j, l})_{l=1}^{M_j}$, 
we can establish the lower bound 
\begin{equation}
    \underline{\eta}_j(m)
    \coloneqq
    \max\bigl\{
    \tilde{\eta}_j(m),\eta_j(0),
    \max_{1\leq\ell\leq M_j}
    \{b_{j,\ell}+|a_{j,\ell}^\intercal m|\}
    \bigr\},
    \label{eq:eta-cut-approximation}
\end{equation}
with $a_{j,\ell}\in\partial\eta_j(m_{j,\ell})$ and 
$b_{j,\ell}=\eta_j(m_{j,\ell})-a_{j,\ell}^{\intercal}m_{j,\ell}$, and where 
$\tilde{\eta}_j(m)$ has been included to better 
capture the growth of $\eta_j$.
Note that $\underline{\eta}_j$ can be interpreted as a 
Benders-type lower model of $\eta_j$ \cite{Diniz+2019}.

\begin{proposition}[Tighter approximation]
    \label{prop:eta-cut-approximation}
    The function $\underline{\eta}_j$ is convex, even, and globally
    Lipschitz. Moreover, for every $m\in\R^{n_x}$,
    \begin{equation}
        0\leq \eta_j(m)-\underline{\eta}_j(m)
        \leq \eta_j(m)-\tilde{\eta}_j(m)
        \leq \eta_j(0).
        \label{eq:eta-cut-uniform-bound}
    \end{equation}
    In particular, $\underline{\eta}_j$ is a global lower bound for
    $\eta_j$ and is pointwise at least as tight as $\tilde{\eta}_j$. It is exact
    at $m=0$ and at every $\pm m_{j,\ell}$.
\end{proposition}

\begin{proof}
    The subgradient inequality and evenness of $\eta_j$ give
    $b_{j,\ell}+|a_{j,\ell}^\intercal m|\leq\eta_j(m)$.
    Also, Proposition~\ref{prop:psi_approx} gives
    $\tilde\eta_j(m)\leq\eta_j(m)$, while convexity and evenness give
    $\eta_j(0)\leq\eta_j(m)$. Thus
    $\tilde\eta_j\leq\underline\eta_j\leq\eta_j$; the remaining bound
    in \eqref{eq:eta-cut-uniform-bound} follows from
    Proposition~\ref{prop:psi_approx} and
    $\eta_j(0)=\rho(\|w\|_{P_j}^2)$.
    Each cut is convex, even, and globally Lipschitz, as is the finite
    seminorm $\tilde\eta_j$ on $\R^{n_x}$; their finite maximum with
    the constant $\eta_j(0)$ has the same properties.
    The constant term makes the bound exact at $0$. At $m_{j,\ell}$,
    $a_{j,\ell}^\intercal m_{j,\ell}+b_{j,\ell}=\eta_j(m_{j,\ell})$;
    evenness gives exactness also at $-m_{j,\ell}$.
\end{proof}

By applying the same arguments as in 
Proposition \ref{prop:sum-seminorms-desiderata}
and using Proposition \ref{prop:eta-cut-approximation}
we can prove the following result.

\begin{proposition}
    \label{prop:hyperplane-method-d1d2}
    Let us define $s_{0}^{*}(x) = 0$ and, recursively,
    \begin{equation}
        s_{j+1}^{*}(x)
        =s_j^{*}(\bar{A}_{j+1}x)
        +\underline{\eta}_j(\bar{A}_{j+1}x),
        \label{eq:s-cut}
    \end{equation}
    and $\widehat{V}_j^{*}(x) = \|x\|_{P_j}^2 + s_j^{*}(x)$,
    where $\underline{\eta}$ is given in 
    \eqref{eq:eta-cut-approximation}. 
    Let $\widehat\kappa_{j+1}^{*}$ be the corresponding 
    controllers, that is 
    \begin{align}
        \kappa^*_{j+1}(x)
        {}={}
        K_{j+1}x {}+{}
        \hspace{-2.5em}
        \argmin_{
            \substack{
                v\in\R^{n_u} 
                \\ 
                z_{i} = \bar{A}_{i+1}z_{i+1}, 
                i=j-1,\ldots, 0
                \\
                z_j = \bar{A}_{j+1}x + Bv
            }
            }
        \hspace{-2.5em}
        \|v\|_{\bar{R}_j}^2
        + \sum_{i=0}^{j} \underline{\eta}_i(z_i).
        \label{eq:kappa_ast}
    \end{align}
    Then, for every $j$,
    $\widehat V_j^{*}$ and $\widehat\kappa_{j+1}^{*}$
    satisfy desiderata \hyperlink{D1}{D1} and \hyperlink{D2}{D2}, with
    $\widehat{\pi}_{j+1}=\widehat\kappa_{j+1}^{*}$.
\end{proposition}

\begin{remark}[Computational complexity]
    In the Gaussian case with the AV@R risk, 
    the minimization in \eqref{eq:kappa_ast} can be cast as 
    an \ac{SOCP} with $O(\sum_i M_i)$ linear inequalities.    
\end{remark}

To obtain subgradients of $\eta_j$ we turn our attention 
to the dual representation of $\rho$. 
Suppose $\rho$ has dual set $\A$. Then, 
\begin{subequations}    
    \label{eq:eta-risk-envelope}
    \begin{equation}
        \eta_j(m) 
        {}={} 
        \sup_{(a, b)\in\mathcal{C}_j}
        \{a^\intercal m + b\},
    \end{equation}
    where $\mathcal{C}_j$ is the convex set 
    \begin{equation}
        \mathcal{C}_j 
        {}={}
        \left\{
            (2P_j \E[yw], \E[y\|w\|_{P_j}^2]),
            y \in \A
            \right\}.
        \end{equation}
\end{subequations}
The set $\A$ is infinite-dimensional,
but $\mathcal{C}_j \subseteq \R^{n_x + 1}$.
For any $(a, b) \in \mathcal{C}_j$, 
$a^\intercal m + b$ is an affine minorant of $\eta_j$.
Suppose that given a point $\hat{m}$ the supremum 
is attained at a $y^\star$, i.e., 
$\eta_j(\hat{m}) = a_j(y^\star)\hat{m} + b_j(y^\star)$,
where $a_j(y^\star) = 2P_j \E[y^\star w]$, and 
$b_j(y^\star) = \E[y^\star \|w\|_{P_j}^2]$. Then, 
$a_j(y^\star)m + b_j(y^\star)$ is a supporting hyperplane, 
and, by definition, $a_j(y^\star) \in \partial \eta_j(\hat{m})$.
Given $\hat{m} \in \R^{n_x}$,
if $\xi_j(\hat{m},w)\in\mathcal{Z}$ and the dual set 
$\A\subseteq\mathcal Z^*$ is w$*$-compact, 
then the dual supremum at $\hat{m}$ is attained.
We will now look at two specific cases of the risk measures 
$\AVaR_\alpha$ and $\EVaR_\alpha$.

\textbf{Average value-at-risk:} 
The average value-at-risk is finite on 
$\mathcal{Z}=\mathcal{L}^1$ with dual set
\begin{equation}
    \label{eq:avar_dual_set}
    \A_\alpha
    {}={}
    \left\{y\in \mathcal{L}^{\infty}:
    0 \leq y \leq \tfrac{1}{1-\alpha}\ \text{a.s.},\
    \E[y] = 1\right\},
\end{equation}
which is w$*$-closed \cite[Sect. 6.3.]{ShapiroDentcheva2021}
and a subset of $\{y\in\mathcal{L}^{\infty}: \|y\|_\infty \leq \tfrac{1}{1-\alpha}\}$, which is w$*$-compact by the 
Banach-Alaoglu theorem, so $\A_\alpha$ is w$*$-compact,
therefore, the maximum is attained by a $y^\star$.
Indeed, given an anchor point $\hat{m}\in\R^{n_x}$ for which
$\xi_j(\hat m,w)$ has a continuous distribution,
by \cite[Eq.~(6.81)]{ShapiroDentcheva2021} the supremum is 
attained at 
\begin{equation}
    y_j^\star
    =\tfrac{1}{1-\alpha}
    1_{\{
      \xi_j(\hat{m},w)\geq t_j^\star
    \}},\;
    \text{a.s.}
    \label{eq:y_star_avar}
\end{equation}
where $t_j^\star=\inf \{t: \prob[\xi_j(\hat{m},w) \leq t] \geq \alpha\}$ is the $\alpha$-quantile of $\xi_j(\hat{m},w)$.
If $w$ is Gaussian, $t_j^\star$ is the quantile 
of a generalized $\chi^2$ distribution \cite{Das2025} 
and can be computed efficiently. 
Lastly, given $y_j^\star$ in \eqref{eq:y_star_avar},
the coefficients can be written as tail expectations:
\begin{subequations}
    \label{eq:ab_tail_expectations}
    \begin{align}
        a_j(y_j^\star)
        {}={}&
        2P_j\E[w{}\mid{}\xi_j(\hat{m},w)\geq t_j^\star],
        \\
        b_j(y_j^\star)
        {}={}&
        \E[\|w\|_{P_j}^2 {}\mid{}\xi_j(\hat{m},w)\geq t_j^\star].
    \end{align}
\end{subequations}
These tail expectations can be computed using Monte Carlo or quadrature
methods \cite{ArismendiBroda2017}, providing the coefficients needed to construct
supporting hyperplanes of $\eta_j$ at selected anchor points.
% This construction is somewhat reminiscent of the cut-based 
% lower approximations used in \ac{SDDP} \cite{GuiguesRomisch2012}.

\textbf{Entropic value-at-risk:} 
Following \cite[Thm.~4.9]{JavidPichler2017}, 
$\EVaR_\alpha$ with $\alpha\in[0, 1)$ 
is defined on the exponential Orlicz space $\mathcal{Z} = \mathcal{L}^\Phi {}\coloneqq{} \{y: \mathcal{F}\text{-measurable}, \exists t > 0: \E e^{t|y|}<\infty\}$, rather than 
on an $\mathcal{L}^p$ space, and its dual set is 
\begin{equation}
    \A = \{y\in\mathcal{Z}^*: y\geq 0 \text{ a.s.}, \E[y]=1, \E[y\log y] \leq c_\alpha\},
\end{equation}
where $c_\alpha = -\log(1-\alpha)$ and $\mathcal{Z}^*$ is 
the complementary Orlicz space $\mathcal{Z}^* {}\coloneqq{} \{y: \mathcal{F}\text{-measurable}, \E [|y| \log(1+|y|)]<\infty\}$.
Let $\Lambda_Z(t) = \log \E e^{tZ}$. If the minimum is 
attained in \eqref{eq:evar-definition} at a point $t^\star > 0$,
then by the dominated convergence theorem $\Lambda_Z$ is 
differentiable at $t^\star$ and the optimality condition is 
\begin{equation}
    t^\star \Lambda_Z'(t^\star) - \Lambda_Z(t^\star) = c_\alpha.
\end{equation}
Given that $t$ is a scalar, this equation can be solved by 
bijection. Let $y^\star = \exp(t^\star Z - \Lambda_Z(t^\star))$;
then $\E[y^\star]=1$, $y^\star \geq 0$, and 
$\E[y^\star \log y^\star] = c_\alpha$, so $y^\star \in \A$, 
and $\EVaR_\alpha[Z] = \E[y^\star Z]$.

The resulting optimal dual variable $y^\star$ can then
be used to construct a supporting hyperplane of $\eta_j$
at the selected anchor point $\hat{m}$.

\begin{remark}[Selection of anchor points]\label{remark:anchor-point-selection}
    The choice of the anchor points in \eqref{eq:eta-cut-approximation}
    can significantly affect the quality of the approximation.
    Here we propose an iterative refinement method.
    At refinement round $r$ we have a set of candidate points 
    $\mathcal{C}_j^{[r]}$ and a set of selected points $\mathcal{M}_j^{[r]}$.
    We initialize $\mathcal{C}_j^{[0]}=\mathcal{M}_j^{[0]}=\varnothing$,
    and $\widehat{\kappa}_{j+1}^{[0]}=\widehat\kappa_{j+1}$ from 
    \eqref{eq:kappa_hat_jp1}. Set
    $\underline\eta_j^{(0)}(m)=\max\{\tilde\eta_j(m),\eta_j(0)\}$.

    At every round $r\geq1$ we generate states $x_{j,p}^{[r]}$ by
    simulating \eqref{eq:system-dynamics} with the approximate
    controllers from round $r-1$ and initial states
    sampled from a compact set $\mathcal{X}\subseteq\R^{n_x}$.
    The corresponding nominal successor states are
    \(
        c_{j,p}^{[r]}
        =Ax_{j,p}^{[r]}
        +B\widehat\kappa_{j+1}^{[r-1]}(x_{j,p}^{[r]}).
    \)
    These are inserted into the set of candidate points,
    $\mathcal{C}_j^{[r]} =\mathcal{C}_j^{[r-1]} \cup\{\pm c_{j,p}^{[r]}\}_p$.
    We carry forward the selected points by setting
    $\mathcal{M}_j^{[r]}=\mathcal{M}_j^{[r-1]}$.

    At each round we select at most $M_j$
    new anchor pairs
    from $\mathcal{C}_j^{[r]}$ using the
    gap $\delta_{j,\ell}(m) = \eta_j(m)-\underline\eta_j^{(\ell-1)}(m)$,
    where $\ell$ counts selected pairs across all rounds. Provided
    unselected candidates remain, their largest gap exceeds
    $\epsilon_{\rm gap}$, and neither the batch nor the total budget
    $M_j$ has been reached,
    we choose
    \begin{equation}
        m_{j,\ell}
        \in\operatorname*{arg\,max}_{m\in\mathcal{C}_j^{[r]}\setminus
        \mathcal{M}_j^{[r]}}
        \delta_{j,\ell}(m),
        \label{eq:greedy-anchor-selection}
    \end{equation}
    and insert $\pm m_{j, \ell}$ into $\mathcal{M}_j^{[r]}$.
    After choosing $a_{j,\ell}\in\partial\eta_j(m_{j,\ell})$ and setting
    $b_{j,\ell}=\eta_j(m_{j,\ell})-a_{j,\ell}^{\intercal}m_{j,\ell}$,
    we update
    \begin{equation}
        \underline{\eta}_j^{(\ell)}(m)
        =\max\bigl\{\underline{\eta}_j^{(\ell-1)}(m),
        b_{j,\ell}+|a_{j,\ell}^{\intercal}m|\bigr\}.
    \end{equation}
    We repeat this step while these conditions hold. If any cuts were
    added, we recompute the approximate controllers
    $\widehat\kappa_{j+1}^{[r]}$. We stop when a round adds no cuts or
    the total budget $M_j$ is reached, otherwise we perform another round.
\end{remark}

Lastly, in Section \ref{section:value-funct-approx}
we established an $O(\|x\|)$-bound for the 
centered error as in \eqref{eq:centered-error-Ox-bound}
when $w$ is Gaussian and $\rho$ satisfies 
Assumption \ref{assm:rho_disturbance}(i),
which guarantees the 
\textit{dilatation monotonicity} condition 
\cite[Prop. 5.6]{Bellini+2021}, namely 
\begin{equation}
    \rho(\E[Z {}\mid{} \mathcal{G}])
    \leq
    \rho(Z),
    \label{eq:dilatation-monotonicity}
\end{equation}
for every $Z\in\mathcal{Z}$ and sub-$\sigma$-algebra 
$\mathcal{G} \subseteq \mathcal{F}$ such that $\E[Z {}\mid{} \mathcal{G}] \in \mathcal{Z}$. 
Both $\AVaR_\alpha$ and $\EVaR_\alpha$, for $\alpha\in[0, 1)$ 
satisfy this property on their respective domains.
% that is 
% $\mathcal{Z} = \mathcal{L}^1$ for $\AVaR_\alpha$ and the 
% exponential Orlicz space for $\EVaR_\alpha$.

\begin{proposition}[Quadratic centered error]%
    \label{prop:quadratic-centered-error}
    Let $\rho$ be a risk measure satisfying 
    Assumption \ref{assm:rho_disturbance}(i) 
    and \eqref{eq:dilatation-monotonicity}.
    Suppose $w\sim \mathcal{N}(0, \Sigma_w)$ with $\Sigma_w \succ 0$. 
    Let $e_j$ be the centered error that corresponds 
    to the hyperplane-based approximate value 
    function $\widehat{V}_j^*$ of 
    Proposition \ref{prop:hyperplane-method-d1d2}.
    Then, for every $j$ there is $c_j \geq 0$ such that for all $x$,
    \begin{equation}
        |e_j(x)| \leq c_j \|x\|^2.
        \label{eq:quadratic-centered-error}
    \end{equation}
\end{proposition}
\begin{proof}
    The proof can be found in the \hyperlink{link:appendix:proof-quadratic-upper-bound}{Appendix}.
\end{proof}

\subsection{Elliptically symmetric disturbances}
The approach of Section \ref{sec:finer-approximation}
exploits the symmetry of $w$ 
(Assumption \ref{assm:rho_disturbance}) 
and approximates $\eta_j$ at points $\pm\hat{m}$.
This can be seen as an approximation of the \textit{tightest 
conic approximation} given in the following theorem.

\begin{theorem}[Tightest conic approximation]
    \label{thm:tightest-conic-approximation}
    Let $a^\intercal m+b$ be an affine minorant of $\eta_j$,
    let $\eta_j^*$ be the convex conjugate of $\eta_j$,
    and define the sub-level set of $\eta_j^*$
    \begin{equation}
        K_b = \{p\in\R^{n_x} : \eta_j^*(p) \leq -b\}.
        \label{eq:Kb_definition}
    \end{equation}
    Define the supporting conic model 
    \begin{equation}
        \eta_j^{\rm con}(m) {}={} b + \sigma_{K_b}(m),
    \end{equation}
    where $\sigma_{K_b}$ is the support function of $K_b$.
    Then, (i) $\eta_j^{\rm con}(m) \leq \eta_j(m)$ for all $m$,
    and (ii) for every finite, convex, positively homogeneous
    $\phi:\R^{n_x}\to\R$ such that $b+\phi\leq\eta_j$, one has
    $b+\phi\leq\eta_j^{\rm con}$. Thus $\eta_j^{\rm con}$ is the
    largest convex conical lower model with fixed intercept $b$.
    If $a^\intercal\hat m+b=\eta_j(\hat m)$ for some $\hat m$,
    then $\eta_j^{\rm con}(\hat m)=\eta_j(\hat m)$.
\end{theorem}

\begin{proof}
    It is $a^\intercal m + b \leq \eta_j(m)
    \Leftrightarrow a^\intercal m -\eta_j(m) \leq -b$, and taking a supremum over $m$
    gives $\eta_j^*(a)\leq-b$,
    so $a\in K_b$. 
    Conversely, for $p\in K_b$, by the definition 
    of the convex conjugate
    \begin{equation}
        p^\intercal m-\eta_j(m)
        {}\leq{}
        \eta_j^*(p)
        {}\leq{}
        -b.
    \end{equation}
    Taking the supremum over $p\in K_b$ proves (i).
    If the affine minorant touches at $\hat m$, equality follows
    because $a\in K_b$. Moreover,
    $a^\intercal m\leq\sigma_{K_b}(m)\leq\eta_j(m)-b$,
    so $\sigma_{K_b}$ is finite, convex, and positively homogeneous.

    Let $b+\phi$ be as in (ii). Since $\phi$ is finite, convex, and
    positively homogeneous, it is a closed sublinear function.
    By \cite[Sec.~3.1, Thm.~3.1.1]{HiriartUrrutyLemarechal2001},
    $\phi=\sigma_C$ for some nonempty closed convex set $C$.
    For every $p\in C$, one has
    $b+p^\intercal m\leq b+\phi(m)\leq\eta_j(m)$ for all $m$, so 
    $\eta_j(m) - p^\intercal m \geq -b $ 
    and by taking a supremum over $m$ 
    we get $p\in K_b$, therefore $C\subseteq K_b$ and
    $\phi=\sigma_C\leq\sigma_{K_b}$.
\end{proof}

Indeed, in the approximation of Section \ref{sec:finer-approximation}
each pair of supporting hyperplanes corresponds to taking an 
inner approximation of the form $\conv \{\pm a\} \subseteq K_b$.
We will now see that under the assumption that $w$ is elliptically 
symmetric and some regularity assumptions, we can obtain a less 
conservative inner approximation of $K_b$.

\begin{assumption}\label{assm:elliptic_symmetry}
    The disturbance $w$ is elliptically symmetric, i.e., the \ac{pdf} of $w$ 
    can be written as $p_w(w) = h(\|w\|_S)$ for some $S\succ 0$.
\end{assumption}

For notational convenience, for vector 
$x, y\in \R^n$ let us define
$[\![x, y]\!] = \sum_i |x_i y_i|$.

\begin{proposition}\label{prop:elliptic-symmetry-cut}
    Suppose Assumption \ref{assm:elliptic_symmetry} holds and let $S^{-\nicefrac{1}{2}}P_jS^{-\nicefrac{1}{2}} = U\Lambda U^\intercal$,
    where $U$ is orthogonal and $\Lambda$ is diagonal, and let
    $T = U^\intercal S^{\nicefrac{1}{2}}$.
    Let $a^{\intercal}m + b$ be an affine minorant of 
    $\eta_j$, i.e., $a^{\intercal}m + b \leq \eta_j(m)$
    for all $m$. Then,
    \begin{equation}
        \eta_j(m) 
        {}\geq{} 
        b 
        {}+{} 
        \llbracket 
            U^\intercal S^{-\nicefrac{1}{2}}a, Tm
        \rrbracket.
        \label{eq:cool_tighter_bound}
    \end{equation}
    If at a point $\hat{m}$ we have 
    $a^{\intercal}\hat{m} + b {}={} \eta_j(\hat{m})$,
    then the bound \eqref{eq:cool_tighter_bound} is exact 
    at every point $T^{-1}DT\hat{m}$
    where $D \in \mathfrak{D} {}\coloneqq{}
    \{\diag(\pm 1, \ldots, \pm 1)\}$.
\end{proposition}

\begin{proof}
    Consider the change of variables $z=Tw$ and note 
    that $\|w\|_S = \|z\|$ and the distribution of $z$
    becomes 
    \begin{subequations}
        \begin{align}
            p_z(z) 
            {}={}& 
            |\det T^{-1}|p_w(T^{-1}z)
            \\
            {}={}&
            (\det S)^{-\nicefrac{1}{2}}h(\|z\|).
        \end{align}
    \end{subequations}
    Let $D \in \mathfrak{D}$.
    Given the density of $z$ we see that 
    $z\overset{\rm d}{\sim}Dz$.
    
    Let $L_D = T^{-1}DT$.
    Since $P_j=T^\intercal\Lambda T$, we have
    $\xi_j(m)=z^\intercal\Lambda z+2(Tm)^\intercal\Lambda z$. 
    As $D\Lambda=\Lambda D$,
    law invariance of $\rho$ gives
    $\eta_j(L_Dm)=\eta_j(m)$.

    Let $a_D=L_D^\intercal a$.
    The assumed minorant and this symmetry imply
    $b+a_D^\intercal m=b+a^\intercal L_Dm\leq\eta_j(m)$;
    hence, by virtue of \eqref{eq:Kb_definition}, $a_D\in K_b$
    for all $D$.
    With $q=T^{-\intercal}a=U^\intercal S^{-\nicefrac{1}{2}}a$,
    Theorem~\ref{thm:tightest-conic-approximation}(i) gives
    \begin{align}
        \eta_j(m)
        {}\geq{}& 
        b+\sigma_{K_b}(m)
        {}\geq{} 
        b+\max_D a_D^\intercal m 
        \notag
        \\
        {}={}& 
        b+\max_D q^\intercal DTm
        {}={} 
        b+\llbracket q,Tm\rrbracket.
    \end{align}
    This proves \eqref{eq:cool_tighter_bound}.
    Suppose $a^{\intercal}\hat{m} + b {}={} \eta_j(\hat{m})$ 
    at a point $\hat{m}$.

    If $a^\intercal\hat m+b=\eta_j(\hat m)$, then the established
    lower bound and the choice $D=I$ give
    \begin{align}
        \eta_j(\hat m)
        {}\geq{}& 
        b + \llbracket q,T\hat m\rrbracket
        \\
        {}\geq{}& 
        b+q^\intercal T\hat m
        =b+a^\intercal\hat m
        =\eta_j(\hat m),
    \end{align}
    so the bound is exact at $\hat m$. For any
    $D \in \mathfrak{D}$, 
    $TL_D\hat m=DT\hat m$, and hence
    $\llbracket q,TL_D\hat m\rrbracket
    =\llbracket q,T\hat m\rrbracket$.
    Since $\eta_j(L_D\hat m)=\eta_j(\hat m)$, the bound is also exact
    at $L_D\hat m$.
\end{proof}

Analogously to \eqref{eq:eta-cut-approximation}, 
Proposition \ref{prop:elliptic-symmetry-cut} produces the
following lower approximation for $\eta_j$
\begin{multline}
    \underline{\eta}_j^{\mathrm{cone}}(m)
    {}={}
    \max
    \Bigl\{
        \tilde{\eta}_j(m),\ \eta_j(0),
    \\
        \max_{1\leq\ell\leq M_j}
        \bigl[b_{j,\ell}
        +\llbracket U^\intercal S^{-1/2}a_{j,\ell},Tm\rrbracket\bigr]
    \Bigr\}.    
\end{multline}

When the disturbance is normally distributed
and $\rho=\AVaR_\alpha$
we can trade the computation of tail expectations
in \eqref{eq:ab_tail_expectations} with non-tight lower bounds
as shown in the following result.

\begin{corollary}[Gaussian--AVaR case]
\label{cor:gaussian-nontight}
    Suppose that $w\sim\mathcal N(0,\Sigma_w)$, $\Sigma_w\succ0$,
    and $\rho=\AVaR_\alpha$ with $\alpha\in[0,1)$.
    Consider the diagonalisation 
    \begin{subequations}
        \begin{equation}
            \Sigma_w^{1/2}P_j\Sigma_w^{1/2} {}={} U\Lambda U^\intercal,        
        \end{equation}
        where $U$ is orthogonal and $\Lambda$ is diagonal, and let 
        \begin{equation}
            T{}={} U^\intercal\Sigma_w^{-1/2}.
        \end{equation}
    \end{subequations}
    For $\alpha\in(0,1)$ and $u$ such that $\|u\|_{\Sigma_w}>0$, let
    $\breve a_j(u)$ and $\breve b_j(u)$  be 
    \begin{subequations}
        \begin{align}
            \breve{a}_j(u) 
            {}={}& 
            \frac{2\kappa_\alpha P_j \Sigma_w u}{\|u\|_{\Sigma_w}},
            \\
            \breve{b}_j(u) 
            {}={}&  
            \trace(P_j \Sigma_w) + \Phi^{-1}(\alpha)\kappa_\alpha
            \frac{\|u\|_{\Sigma_w P_j \Sigma_w}^2}{\|u\|_{\Sigma_w}^2}.
        \end{align}
    \end{subequations}
    For $\alpha=0$, take $\breve a_j(u)=0$ and
    $\breve b_j(u)=\trace(P_j\Sigma_w)$.
    Then, for every $m\in\R^{n_x}$,
    \begin{subequations}
        \begin{align}
            \eta_j(m)
            {}\geq{}&
            \breve b_j(u)+
            \llbracket
                U^\intercal\Sigma_w^{1/2}\breve a_j(u),  Tm
            \rrbracket
            \\
            {}\geq{}& \breve b_j(u)+|\breve a_j(u)^\intercal m|.
            \label{eq:gaussian-conic-cut}
        \end{align}
    \end{subequations}
\end{corollary}

% \begin{corollary}[Special case]
%     \label{cor:gaussian-nontight}
%     Suppose that $w\sim \mathcal{N}(0, \Sigma_w)$ and $\rho=\AVaR_\alpha$ with $\alpha\in[0, 1)$.
%     Then, for every
%     $m\in\R^{n_x}$ and $u\in\R^{n_x}$ such that
%     $\|u\|_{\Sigma_w}>0$,
%     $\eta_j(m) \geq \breve{a}_j(u)^\intercal m + \breve{b}_j(u)$,
%     where, for $\alpha\in(0,1)$,
%     \begin{subequations}
%         \begin{align}
%             \breve{a}_j(u) 
%             {}={}& 
%             \frac{2\kappa_\alpha P_j \Sigma_w u}{\|u\|_{\Sigma_w}},
%             \\
%             \breve{b}_j(u) 
%             {}={}&  
%             \trace(P_j \Sigma_w) + \Phi^{-1}(\alpha)\kappa_\alpha
%             \frac{\|u\|_{\Sigma_w P_j \Sigma_w}^2}{\|u\|_{\Sigma_w}^2}.
%         \end{align}
%     \end{subequations}
%     For $\alpha=0$, take $\breve a_j(u)=0$ and
%     $\breve b_j(u)=\trace(P_j\Sigma_w)$.
% \end{corollary}
%
\begin{proof}
    If $\alpha=0$, then $\AVaR_0=\E$ and
    $\eta_j(m)=\trace(P_j\Sigma_w)$, so the assertion holds.
    Suppose $\alpha\in(0,1)$ and let $\A$ be the dual set of $\AVaR_\alpha$.
    Let $n=u^\intercal w/\|u\|_{\Sigma_w}\sim\mathcal N(0,1)$, and
    consider the event $E=\{n\geq\Phi^{-1}(\alpha)\}$. 
    Since $\prob(E)=1-\alpha$,
    \begin{equation}
        \breve{y}=\tfrac{1}{1-\alpha}1_{E}
        {}\in{} \A.
        \label{eq:gaussian-cut-density}
    \end{equation}
    By \eqref{eq:eta-risk-envelope}, this gives
    $\eta_j(m)
        \geq 2m^\intercal P_j\E[\breve{y}w]
        +\E[\breve{y}\|w\|_{P_j}^2]$.
    We now need to determine $\E[\breve{y}w]$ and 
    $\E[\breve{y}\|w\|_{P_j}^2]$.

    From \eqref{eq:avar-tail-definition}, $\E[n {}\mid{} E] = \AVaR_\alpha[n] = \kappa_\alpha$.
    Note that $(w, n)$ is jointly normal
    with $\cov(w, n) = \Sigma_wu / \|u\|_{\Sigma_w}$, so
    $\E[w {}\mid{} n] = \frac{\Sigma_w u}{\|u\|_{\Sigma_w}}n$.
    By the tower property of total expectation,
    \begin{equation}
        \E[w {}\mid{} E] = \E[\E[w {}\mid{} n]{}\mid{}E]
        =\frac{\kappa_\alpha\Sigma_wu}{\|u\|_{\Sigma_w}}.
    \end{equation}
    From \eqref{eq:gaussian-cut-density}, 
    $\E[\breve{y}w]=\E[w\mid E]$, therefore
    \begin{equation*}
        2P_j\E[\breve{y} w]
        =\frac{2\kappa_\alpha P_j\Sigma_wu}
              {\|u\|_{\Sigma_w}}
        =\breve a_j(u).
    \end{equation*}
    Next, $\E[n^2 {}\mid{} E] = \E[n^2 1_E] / \prob[E]$,
    and since $n \sim \mathcal{N}(0, 1)$ and $\prob[E] = 1-\alpha$,
    \begin{align}
        \E[n^2 {}\mid{} E] 
        {}={} 
        \tfrac{1}{1-\alpha}
        \int_{\Phi^{-1}(\alpha)}^{\infty} 
        \hspace{-0.5em}
        \tau^2 \phi(\tau)\dd \tau
        {}={}
        1 + \kappa_\alpha \Phi^{-1}(\alpha),
        \label{eq:expectation-nsq-given-E}
    \end{align}
    where the last equality follows by integration by parts using
    $\phi'(\tau)=-\tau\phi(\tau)$.
    Again by the tower property $\E[ww^\intercal {}\mid{} E] = \E[ \E[ww^\intercal {}\mid{} n] {}\mid{} E]$, and using the conditional variance we obtain
    \begin{align}
        \E[ww^\intercal {}\mid{} n]
        {}={}&
        \var(w {}\mid{} n) + \E[w{}\mid{}n]\E[w{}\mid{}n]^\intercal
        \notag 
        \\
        {}={}&
        \Sigma_w - \frac{\Sigma_w uu^\intercal \Sigma_w}{\|u\|_{\Sigma_w}^2} 
        + \frac{\Sigma_w uu^\intercal\Sigma_w}
        {\|u\|_{\Sigma_w}^2}n^2.
    \end{align}
    Conditioning on $E$ and using 
    \eqref{eq:expectation-nsq-given-E} we obtain 
    \begin{equation}
        \E[ww^\intercal {}\mid{} E] {}={} \Sigma_w + \tfrac{1}{\|u\|_{\Sigma_w}^2}\Phi^{-1}(\alpha)\kappa_\alpha \Sigma_w uu^\intercal \Sigma_w.
    \end{equation}
    Hence $\E[\breve{y} \|w\|_{P_j}^2]
    = \trace(P_j \E[ww^\intercal {}\mid{} E]) = \breve{b}_j(u)$.
    Thus $\breve a_j(u)^\intercal m+\breve b_j(u)$ is an affine
    minorant of $\eta_j$.
    A Gaussian density satisfies Assumption~\ref{assm:elliptic_symmetry}
    with $S=\Sigma_w^{-1}$, so Proposition~\ref{prop:elliptic-symmetry-cut}
    gives the first inequality. The second follows from
    \begin{align}
        |\breve a_j(u)^\intercal m|
        {}={}&
        \bigl|\bigl(U^\intercal\Sigma_w^{1/2}\breve a_j(u)\bigr)^\intercal Tm\bigr|
        \\
        {}\leq{}&
        \llbracket U^\intercal\Sigma_w^{1/2}\breve a_j(u),Tm\rrbracket.
    \end{align}
    This completes the proof.
\end{proof}

For the Gaussian--$\AVaR$ case, Corollary \ref{cor:gaussian-nontight} allows for
the closed-form construction of
$\breve{a}_{j, l} = \breve{a}_{j}(u_l)$ and 
$\breve{b}_{j, l} = \breve{b}_{j}(u_l)$ for a collection of points $u_l$.
% \cmnt{The points $u_l$ can be random samples from the unit sphere or adaptively as in Remark \ref{remark:anchor-point-selection}.}

\section{Structured approximate dynamic programming}\label{sec:structured-adp}
In the previous section we presented approximation methods 
for $\eta_j$ while omitting the term $r_j$ in $\psi_j$ in \eqref{eq:psij}.
We now develop a structure-preserving \ac{aDP} approach
that exploits Theorem~\ref{thm:Vjstar_structure}
by retaining the known quadratic component of the
optimal value function and approximating only its residual.
This \emph{grey-box} approximation leverages the analytical
structure of the problem, and as we shall demonstrate in 
Section~\ref{sec:simulations} it can improve accuracy
and generalization, as illustrated.

Specifically, we approximate $V_j^\star$ via \ac{aDP}
\cite{Bertsekas2012_DPvol2_aDP},
% Specifically, 
% since we know from Theorem \ref{thm:Vjstar_structure} 
% that $V_j^\star$ possesses a certain structure, we can solve 
% the problem via \ac{aDP} \cite{Bertsekas2012_DPvol2_aDP}, 
using the candidate value functions of the form 
\begin{equation}
    \widetilde{V}_j(x) = \|x\|_{P_j}^2 + f_j(x; \theta),
\end{equation}
where $f_j(x; \theta)$ is a parametric function, with parameters 
$\theta \in \R^{n_\theta}$---to be learned---designed such that it 
is convex, even, and globally Lipschitz.
This can be an \ac{ICNN}
\cite{AmosXiKotler2017} with symmetrisation to enforce 
evenness. 

% Describe aDP procedure
The \ac{aDP} procedure constructs recursively approximate value 
functions $\widetilde{V}_j$ with $\widetilde{V}_0 = V_f$.
Given test points $(x_i)_i$ \ac{aDP} evaluates 
$y_{j+1, i} {}\approx{} \T \widetilde{V}_j(x_i)$
by approximately solving 
$\T \widetilde{V}_j(x_i) = \inf_{u}\{\ell(x, u) + \rho[\widetilde{V}_j(Ax_i + Bu + w)]\}$.
Then $f_{j+1}(x; \theta)$ is fitted to $y_{j+1, i} - \|x_i\|_{P_{j+1}}^2$.
This approach offers greater flexibility compared to the 
methods above, but the approximate value functions and the 
corresponding controllers cannot be shown to satisfy 
\hyperlink{D1}{D1} and \hyperlink{D2}{D2} globally. 
Lastly, the \ac{aDP} control laws are
\begin{subequations}
    \label{eq:control_law_adp}
    \begin{align}
        \widetilde\kappa_{j+1}(x)
        {}={}&
        K_{j+1}x
        {}+{}
        \argmin_{
            v\in\mathbb R^{n_u}
        }
        \{ \|v\|_{\bar R_j}^2 + \rho [Z_j(v; w)] \},
        \label{eq:control_law_adp:a}
        \intertext{where}
        Z_j(v; w)
        {}={}&
        2m_v^\intercal P_jw
        {}+{} \|w\|_{P_j}^2
        {}+{} f_j(m_v+w;\theta_j),
        \label{eq:control_law_adp:b}
        \\
        m_v {}={}& \bar{A}_{j+1}x + Bv.
        \label{eq:control_law_adp:c}
    \end{align}
\end{subequations}
The optimization problem in \eqref{eq:control_law_adp}
is strongly convex, provided $f_j$ is convex.
Compared to the earlier \ac{SOCP} formulations,
the main challenge is that the risk measure, $\rho$, needs to be evaluated from samples.
As an example, in the $\AVaR$ case we can
draw $N_s$ disturbance samples
$w^{(1)},\ldots,w^{(N_s)}$ and estimate the risk
$\rho [Z_j(v; w)]$ in \eqref{eq:control_law_adp:a} by
\cite[Sec.~7.5.1]{ShapiroDentcheva2021}
\begin{equation}
    \widehat{\AVaR}_{\alpha,N_s}(v)
    =\min_{t\in\R}
        t+\frac{\sum\limits_{i=1}^{N_s}(Z_i(v; w^{(i)})-t)_+}{(1-\alpha)N_s}.
    \label{eq:adp-sample-avar}
\end{equation}
It should be noted that for heavy-tailed distributions 
or values of $\alpha$ close to $1$, a large number of 
samples is needed.
Similarly, $\EVaR_\alpha$ can be approximated from samples 
as in \cite[Sec. 3.2]{JavidTafti2019}.

% \cmnt{What about bound on compact sets?}

Lastly, during training the approximate functions $\widehat{V}_j$ 
from Section \ref{section:value-funct-approx} can be used as a regularizer 
to guide training.

Lastly, let $\varsigma(t)=\log(1+\exp(t))$ denote the softplus function. 
We can define a deep convex-even-Lipschitz network by
\begin{subequations}
    \label{eq:adp-deep-convex-even}
    \begin{align}
        z_1(x)
        {}={}& 
        \varsigma(|U_1x|+d_1),\\
        z_{k+1}(x)
        {}={}& 
        \varsigma\bigl(W_kz_k(x)+|U_{k+1}x|+d_{k+1}\bigr),\\
        r^{\rm dcel}_L(x)
        {}={}&
        c+a^\intercal z_L(x)+|v^\intercal x|,
    \end{align}
\end{subequations}
with $W_k\geq0$ and $a\geq0$, componentwise,
for $k=1,\ldots, L-1$.
Here, absolute values in the hidden layers are componentwise and the $U_k\in\R^{m_k\times n_x}$,
$W_k\in\R^{m_{k+1}\times m_k}$, $d_k\in\R^{m_k}$,
$v\in\R^{n_x}$, $a\in\R^{m_L}$, and $c\in\R$---where 
$m_k$ is the width of layer $k$---are 
parameters to be learned. 
The structure of $r^{\rm dcel}_L$ is inspired by 
\acp{ICNN} \cite{AmosXiKotler2017}, while the absolute values 
enforce evenness. 
Its global Lipschitz constant can be controlled by 
imposing suitable norm constraints on its weights 
following~\cite{Gouk+2021}.

\begin{proposition}
    For every $L\in\N_{\geq 1}$, 
    $r^{\rm dcel}_L$ in 
    \eqref{eq:adp-deep-convex-even} is convex, even, and globally Lipschitz.
\end{proposition}

\begin{proof}
    (\textit{Convexity and evenness}) The softplus function is convex, nondecreasing, and nonexpansive ($1$-Lipschitz); indeed,
    $\mathrm{softplus}'(x) = \frac{e^x}{1+e^x} \in (0, 1)$.
    Each component of $x\mapsto |U_kx|$ is convex and even.
    Since $W_k\geq0$ componentwise, induction shows that every
    component of $z_k$ is convex and even. 
    The condition $a\geq0$
    and the convexity and evenness of $x\mapsto|v^\intercal x|$
    give that $r_L^{\rm dcel}$ is convex and even.

    (\textit{Lipschitzianity}) Componentwise, absolute value and softplus are nonexpansive in the Euclidean norm, i.e., $\bigl\||z|-|z'|\bigr\|_2\leq\|z-z'\|_2$ and $\|\varsigma(z)-\varsigma(z')\|_2\leq\|z-z'\|_2$ for all $z,z'$.
    Let $\lambda_1=\|U_1\|_2$ and recursively define 
    $\lambda_{k+1}=\|W_k\|_2\lambda_k+\|U_{k+1}\|_2$.
    By induction, it follows that for all $x,x'$
    \begin{equation}
        \|z_k(x)-z_k(x')\|\leq\lambda_k\|x-x'\|,
    \end{equation}
    therefore,
    \begin{equation}
        |r_L^{\rm dcel}(x)-r_L^{\rm dcel}(x')|
        {}\leq{}
        (\|a\|\lambda_L+\|v\|) \|x-x'\|,
    \end{equation}
    which proves the Lipschitzianity of $r_L^{\rm dcel}$.
\end{proof}

\section{Simulations}\label{sec:simulations}

We shall now demonstrate the effectiveness of the 
value function approximation methods of Section 
\ref{section:value-funct-approx}.
All simulations were performed on a Apple M4 Pro machine 
with 24GB RAM, in Python 3.13.14.

\subsection{System with Gaussian disturbance}\label{sec:simulations:gaussian_dist}
\textbf{System dynamics.} Consider a dynamical system with 
\begin{equation}
    A = \begin{bmatrix}
        1 & 0.2 \\ 0 & 1
    \end{bmatrix},
    B = \begin{bmatrix}
        0.02 \\ 0.2
    \end{bmatrix},
\end{equation}
and $w_k \overset{\text{iid}}{\sim} \mathcal{N}(0, \Sigma_w)$, where $\Sigma_w = \left[\begin{smallmatrix}
    0.01 & 0.003 \\ 0.003 & 0.01
\end{smallmatrix}\right]$, $Q = \diag(1, 0.1)$, $R = 0.1$, and $P_f\approx
\left[\begin{smallmatrix}
4.8174&1.5890\\
1.5890&1.4221
\end{smallmatrix}\right]$ is taken to be 
the solution of the discrete-time algebraic 
Riccati equation.

\textbf{Reference method.} Given that the multistage risk-averse problem does not
admit a closed-form solution, we propose an approximate 
method to be used as a reference for comparison. 
Take a compact sets $X\subseteq \R^{n_x}$ and 
$U\subseteq\R^{n_u}$ and a grid of 
points $\mathcal{G}_x\subseteq X$ on which we will approximate 
$r_j$ and $\psi_j$. 
We apply a Gauss-Hermite quadrature to 
the quantity $2m^\intercal P_jw+\|w\|_{P_j}^2+r_j(m+w)$ 
in \eqref{eq:psij} to approximate $\psi_j$ 
on $\mathcal{G}_x$~\cite[Sec.~3.5(v)]{Olver2010}.
The values $r_j(m+w)$ are determined by linear interpolation
of the $r_j$ values on $\mathcal{G}_x$.
To approximate $r_{j+1}$ for each $x\in \mathcal{G}_x$ 
we use a derivative-free method: we evaluate the objective 
in \eqref{eq:r-j-plus-one}
on the points of a grid $\mathcal{G}_u \subseteq U$.
We then use a golden section search to refine the 
solution~\cite[Sec. 3.4]{KochenderferWheeler2019}.
The procedure is repeated over finer grids until 
the difference between successive approximations 
becomes small. 
This way we construct a reference value function 
$V_j^{\mathrm{ref}}$.

\textbf{Value function approximation error.} To compare an approximate value function 
$\widehat{V}_j(x)$ with the reference 
$V_j^{\mathrm{ref}}$ we define
$\Delta \widehat{V}_j(x) = \widehat{V}_j(x) - \widehat{V}_j(0)$ and 
$\Delta V_j^{\mathrm{ref}}(x) = V_j^{\mathrm{ref}}(x) - V_j^{\mathrm{ref}}(0)$
and define the relative maximum error
\begin{equation}\label{eq:value-relative-error}
    \varepsilon_{\infty,j}^{V, \text{rel}} 
    {}={} 
    \frac{
        \max\limits_{x\in\mathcal G_x} 
        \left| \Delta\widehat{V}_j(x)- \Delta V_j^{\mathrm{ref}}(x) \right|
    }{
        \max\limits_{x\in\mathcal G_x} 
        \Delta V_j^{\mathrm{ref}}(x)
    }\cdot 100\%,
\end{equation}
where $\mathcal{G}_x$ is a uniform grid of 
points in $[-2,2]\times[-1.5,1.5]$.
The results are summarized in Table \ref{tab:value-relative-error}.

\begin{table}[htpb]
    \centering
    \caption{Centered relative value-function errors as defined in \eqref{eq:value-relative-error} for nested $\AVaR_\alpha$ and $\EVaR_\alpha$ at different horizon lengths $j$. Lighter shades indicate smaller errors.}
    \label{tab:value-relative-error}
    \scriptsize
    \setlength{\tabcolsep}{3.2pt}
    \renewcommand{\arraystretch}{1.08}
\begin{tabular}{ccccc@{\hspace{8pt}}cc@{\hspace{8pt}}c}
\toprule
&&\multicolumn{3}{c}{$\AVaR_\alpha$} & \multicolumn{2}{c}{$\EVaR_\alpha$} & \\
\cline{3-5}\cline{6-7}\\
$\alpha$
& $j$
& Sec.~\ref{sec:sum-seminorms}
& Eq.~\eqref{eq:ab_tail_expectations}
& Cor.~\ref{cor:gaussian-nontight}
& Sec.~\ref{sec:sum-seminorms}
& Sec.~\ref{sec:finer-approximation}
& LQ \\
\midrule
0.10 & 1  & \heatA{0.1} & \heatA{0.0} & \heatA{0.0} & \heatA{0.1} & \heatA{0.1} & \heatD{1.3} \\
    %  & 2  & \heatB{0.002} & \heatA{0.000} & \heatA{0.000} & \heatD{0.022} \\
     & 5  & \heatC{0.4} & \heatA{0.1} & \heatA{0.1} & \heatC{0.9} & \heatB{0.4} & \heatE{3.4} \\
     & 10 & \heatC{0.7} & \heatA{0.1} & \heatA{0.1} & \heatD{1.7} & \heatC{0.5} & \heatE{3.7} \\
     & 20 & \heatC{0.9} & \heatB{0.2} & \heatA{0.1} & \heatD{2.3} & \heatC{0.8} & \heatE{3.7} \\
\addlinespace[0.3ex]
0.50 & 1  & \heatB{0.2} & \heatA{0.0} & \heatA{0.1} & \heatC{0.4} & \heatA{0.1} & \heatE{5.2} \\
    %  & 2  & \heatC{0.005} & \heatA{0.001} & \heatA{0.001} & \heatE{0.086} \\
     & 5  & \heatD{1.3} & \heatC{0.5} & \heatC{0.5} & \heatD{2.3} & \heatD{1.1} & \heatF{12.9} \\
     & 10 & \heatD{2.3} & \heatC{0.8} & \heatC{0.7} & \heatE{4.4} & \heatD{1.6} & \heatF{13.7} \\
     & 20 & \heatE{3.2} & \heatD{1.1} & \heatC{0.7} & \heatE{5.8} & \heatD{2.6} & \heatF{13.8} \\
\addlinespace[0.3ex]
0.90 & 1  & \heatC{0.5} & \heatB{0.2} & \heatC{0.4} & \heatC{0.7} & \heatB{0.2} & \heatF{10.8} \\
    %  & 2  & \heatD{0.010} & \heatC{0.005} & \heatC{0.006} & \heatF{0.172} \\
     & 5  & \heatD{2.9} & \heatD{1.9} & \heatD{1.9} & \heatE{4.5} & \heatD{2.7} & \heatF{24.3} \\
     & 10 & \heatE{5.6} & \heatD{2.9} & \heatD{2.8} & \heatE{8.4} & \heatE{3.9} & \heatF{24.9} \\
     & 20 & \heatE{7.4} & \heatE{3.8} & \heatE{3.5} & \heatF{10.6} & \heatE{5.5} & \heatF{25.1} \\
\addlinespace[0.3ex]
0.95 & 1  & \heatC{0.6} & \heatB{0.2} & \heatC{0.4} & \heatC{0.8} & \heatB{0.3} & \heatF{12.4} \\
    %  & 2  & \heatD{0.013} & \heatC{0.007} & \heatC{0.007} & \heatF{0.195} \\
     & 5  & \heatE{3.5} & \heatD{2.5} & \heatD{2.5} & \heatE{5.2} & \heatE{3.3} & \heatF{27.2} \\
     & 10 & \heatE{6.8} & \heatE{3.9} & \heatE{3.8} & \heatE{9.6} & \heatE{4.8} & \heatF{27.7} \\
     & 20 & \heatE{8.9} & \heatE{5.0} & \heatE{4.6} & \heatF{12.1} & \heatE{6.5} & \heatF{27.9} \\
\bottomrule
\end{tabular}
\par\vspace{0.5ex}\hspace{3em}
{\small lower error\enspace
\heatlegendbox{HeatA}\heatlegendbox{HeatB}\heatlegendbox{HeatC}%
\heatlegendbox{HeatD}\heatlegendbox{HeatE}\heatlegendbox{HeatF}%
\enspace higher error}

\end{table}

% \cmnt{\textbf{TODO:} Create a separate table (here) with a comparison 
% between aDP --- structured (DCEL) vs unstructured (MLP). A more 
% thorough comparison of aDP methods will be in a separate 
% subsection.}

% The LQ column gives the risk-neutral Riccati baseline 
% relative to the $\AVaR_\alpha$ reference. 
% The Eq.~\eqref{eq:ab_tail_expectations} column uses 
% Proposition~\ref{prop:elliptic-symmetry-cut} and order-112 Gauss--Hermite 
% quadrature. The two $\AVaR_\alpha$ cut-based columns use $M_j=16$ 
% anchor pairs per stage, selected as in Remark~\ref{remark:anchor-point-selection}; 
% Corollary~\ref{cor:gaussian-nontight} uses one direction per selected anchor, 
% chosen by maximizing its analytic cut over an angular grid. 
% The $\EVaR_\alpha$ Section~\ref{sec:finer-approximation} column 
% retains $M_j=128$ symmetric anchor pairs per stage.

Table~\ref{tab:value-relative-error} reports the centered
relative value-function errors for the different
approximations under $\AVaR_\alpha$ and $\EVaR_\alpha$.
The \ac{LQ} column provides the risk-neutral Riccati baseline
relative to the $\AVaR_\alpha$ reference, while the
sum-of-seminorms approximation is presented in
Section~\ref{sec:sum-seminorms}.
For $\AVaR_\alpha$, the finer approximation based on
Eq.~\eqref{eq:ab_tail_expectations} uses
Proposition~\ref{prop:elliptic-symmetry-cut} and
order-112 Gauss--Hermite quadrature.
Both $\AVaR_\alpha$ cut-based constructions use $M_j=16$
anchor pairs per stage, selected as in
Remark~\ref{remark:anchor-point-selection};
Corollary~\ref{cor:gaussian-nontight} uses one direction
per selected anchor, chosen by maximizing its
analytical cut over an angular grid.
For $\EVaR_\alpha$, the finer approximation of
Section~\ref{sec:finer-approximation} uses $M_j=128$
symmetric anchor pairs per stage.

Table \ref{tab:value-relative-error} shows that the centered 
relative errors generally increase with both the level of 
risk aversion $\alpha$ and the horizon $j$. 
For both $\AVaR$ and $\EVaR$, the finer approximation in 
Section \ref{sec:finer-approximation} outperforms 
the sum-of-seminorms approximation in all cases. 
The two $\AVaR$ cut constructions give similar results, with Corollary \ref{cor:gaussian-nontight} being slightly more 
accurate at longer horizons.

\begin{table}[htpb]
    \caption{Centered relative value-function errors as defined in \eqref{eq:value-relative-error} 
    for nested $\AVaR_\alpha$ at different
    horizon lengths $j$. 
    Lighter shades indicate smaller errors.}
    \centering
    \begingroup
    \let\heatA\undefined
    \let\heatB\undefined
    \let\heatC\undefined
    \let\heatD\undefined
    \let\heatE\undefined
    \let\heatF\undefined
    \let\heatlegendbox\undefined
    % Slide copy of Table I, omitting the Gaussian-only Proposition 9 column.
% AVaR aDP columns: radius 5, 512 training and validation disturbances,
% 1024 training states per stage, seed 20261026.
% Unstructured ICNN: full-value four-layer network; held-out early stopping, at most 4000 epochs for alpha=0.5,0.9; other levels unchanged.
\begingroup
\definecolor{HeatA}{HTML}{FBFDFF}
\definecolor{HeatB}{HTML}{F1F7FC}
\definecolor{HeatC}{HTML}{E6F1F9}
\definecolor{HeatD}{HTML}{D8E9F5}
\definecolor{HeatE}{HTML}{C8DFF0}
\definecolor{HeatF}{HTML}{B5D3EA}
\newcommand{\heatA}[1]{\cellcolor{HeatA}\strut #1}
\newcommand{\heatB}[1]{\cellcolor{HeatB}\strut #1}
\newcommand{\heatC}[1]{\cellcolor{HeatC}\strut #1}
\newcommand{\heatD}[1]{\cellcolor{HeatD}\strut #1}
\newcommand{\heatE}[1]{\cellcolor{HeatE}\strut #1}
\newcommand{\heatF}[1]{\cellcolor{HeatF}\strut #1}
\newcommand{\heatlegendbox}[1]{%
  \begingroup\setlength{\fboxsep}{0pt}%
  \colorbox{#1}{\rule{1.15em}{0pt}\rule{0pt}{1.05ex}}%
  \endgroup}
\scriptsize
\setlength{\tabcolsep}{4pt}
\renewcommand{\arraystretch}{1.0}
\begin{tabular}{cc@{\hspace{6pt}}ccc@{\hspace{6pt}}c}
\toprule
% $\varepsilon_{\infty,j}^{V, \text{rel}}$
&&\multicolumn{3}{c}{$\AVaR_\alpha$}&\\
\cmidrule(lr){3-5}
$\alpha$ & $j$
 & \shortstack{structured\\DCEL$_4$}
 & \shortstack{unstructured\\MLP$_4$}
 & \shortstack{unstructured\\ICNN$_4$} \\
\midrule
\addlinespace[0.25ex]
0.10 & 1  &   \heatA{0.05} & \heatD{1.3} & \heatD{1.8}\\
     & 5  &   \heatA{0.1} & \heatD{2.6} & \heatE{3.4} \\
     & 10 &   \heatB{0.2} & \heatE{4.8} & \heatD{2.9} \\
     & 20 &   \heatB{0.2} & \heatE{5.9} & \heatE{7.3} \\
\addlinespace[0.25ex]
0.50 & 1  &   \heatB{0.3} & \heatD{1.4} & \heatD{1.3} \\
     & 5  &   \heatD{1.2} & \heatD{2.3} & \heatE{3.1} \\
     & 10 &   \heatC{0.9} & \heatE{4.2} & \heatE{4.3} \\
     & 20 &   \heatD{1.0} & \heatE{5.6} & \heatE{6.7} \\
\addlinespace[0.25ex]
0.90 & 1  &   \heatC{0.8} & \heatD{1.2} & \heatD{2.1} \\
     & 5  &   \heatD{2.8} & \heatD{1.9} & \heatE{4.4} \\
     & 10 &   \heatD{2.1} & \heatE{3.2} & \heatE{4.2} \\
     & 20 &   \heatD{1.8} & \heatE{5.7} & \heatE{5.8} \\
\addlinespace[0.25ex]
0.95 & 1  &   \heatC{0.9} & \heatD{1.2} & \heatD{2.6}  \\
     & 5  &   \heatD{2.9} & \heatD{2.2} & \heatE{4.6}  \\
     & 10 &   \heatD{2.9} & \heatE{3.1} & \heatE{3.6}  \\
     & 20 &   \heatD{2.5} & \heatE{4.7} & \heatE{4.8}  \\
\bottomrule
\end{tabular}
\par\vspace{0.5ex}
{\scriptsize lower error\enspace
\heatlegendbox{HeatA}\heatlegendbox{HeatB}\heatlegendbox{HeatC}%
\heatlegendbox{HeatD}\heatlegendbox{HeatE}\heatlegendbox{HeatF}%
\enspace higher error}
\endgroup

    \endgroup
    \label{tab:value-relative-error_adp}
\end{table}

% The proposed structure-preserving
%     DCEL$_4$ approximation is compared against unstructured
%     MLP$_4$ and ICNN$_4$ approximations.

% Table~\ref{tab:value-relative-error_adp} compares the
% proposed structured DCEL approximation in
% \eqref{eq:adp-deep-convex-even} against unstructured
% MLP and ICNN approximations described in
% Section~\ref{sec:structured-adp}.
Table~\ref{tab:value-relative-error_adp} reports the
centered relative value-function errors for the
structured DCEL approximation proposed in
\eqref{eq:adp-deep-convex-even}, compared with
unstructured MLP and ICNN approximations described in
Section~\ref{sec:structured-adp}.
For the numerical comparison, we consider
four-layer architectures, denoted by
DCEL$_4$, MLP$_4$, and ICNN$_4$, respectively.
The structured DCEL approximation generally achieves
lower centered relative errors, with the improvement
becoming more pronounced at longer horizons.
These results illustrate the benefits of exploiting
the analytical structure of the optimal value function
for improving approximation accuracy and generalization.

\begin{figure}
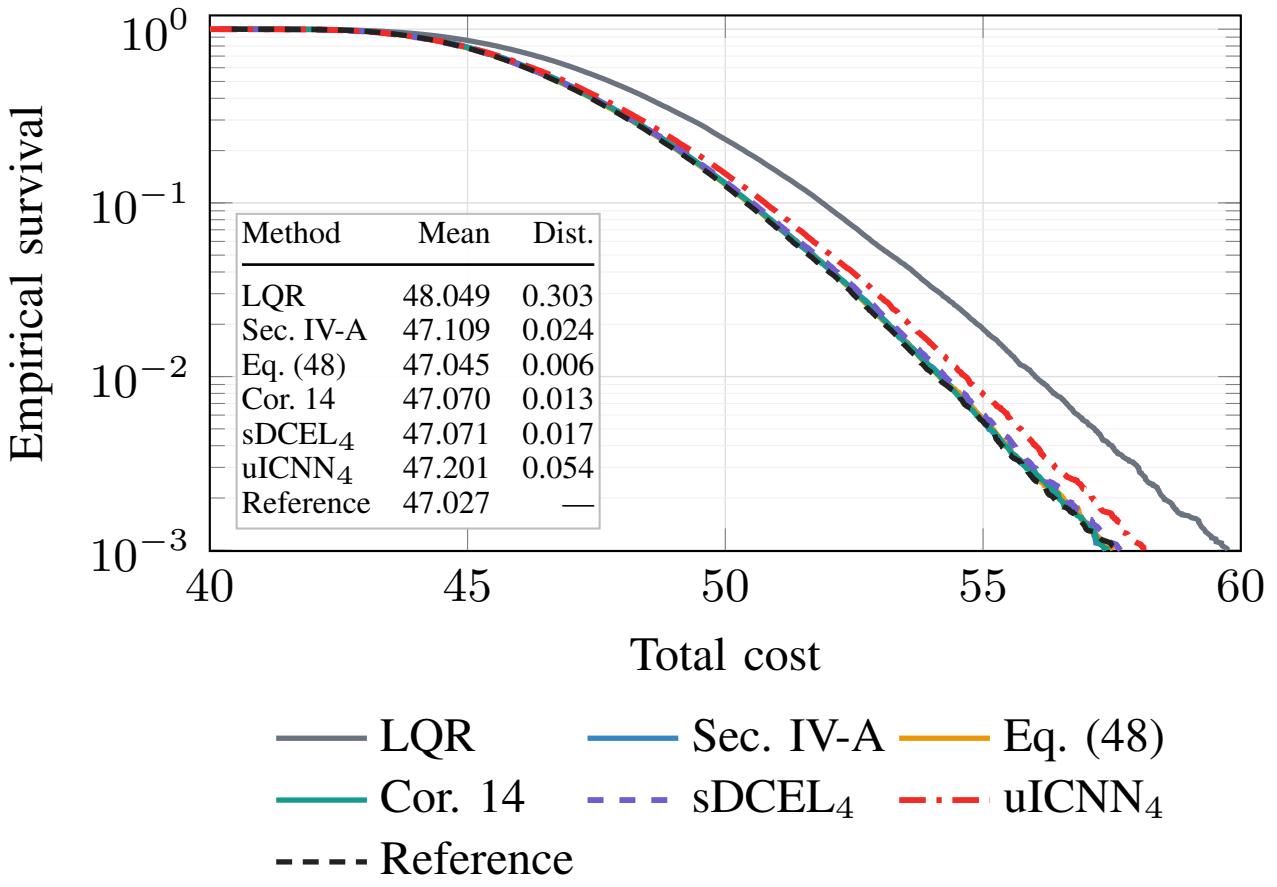

    \paperfigure{adversarial_cost_survival_h15_x2p25_1p5_a05}
    \caption{Empirical total-cost survival functions under each policy's own discrete AVaR adversary, starting from $x_0=(2.25,1.5)^\intercal$ ($j=15$, $\alpha=0.5$, and $50{,}000$ paths per policy). The aDP controllers use $\|x\|_{P_j}^2+r^{\rm dcel}_4(x)$ (sDCEL$_4$) and $g^{\rm icnn}_4(x)$ (uICNN$_4$).}
    \label{fig:adversary-cost-survival-alpha05}
\end{figure}

\begin{figure}
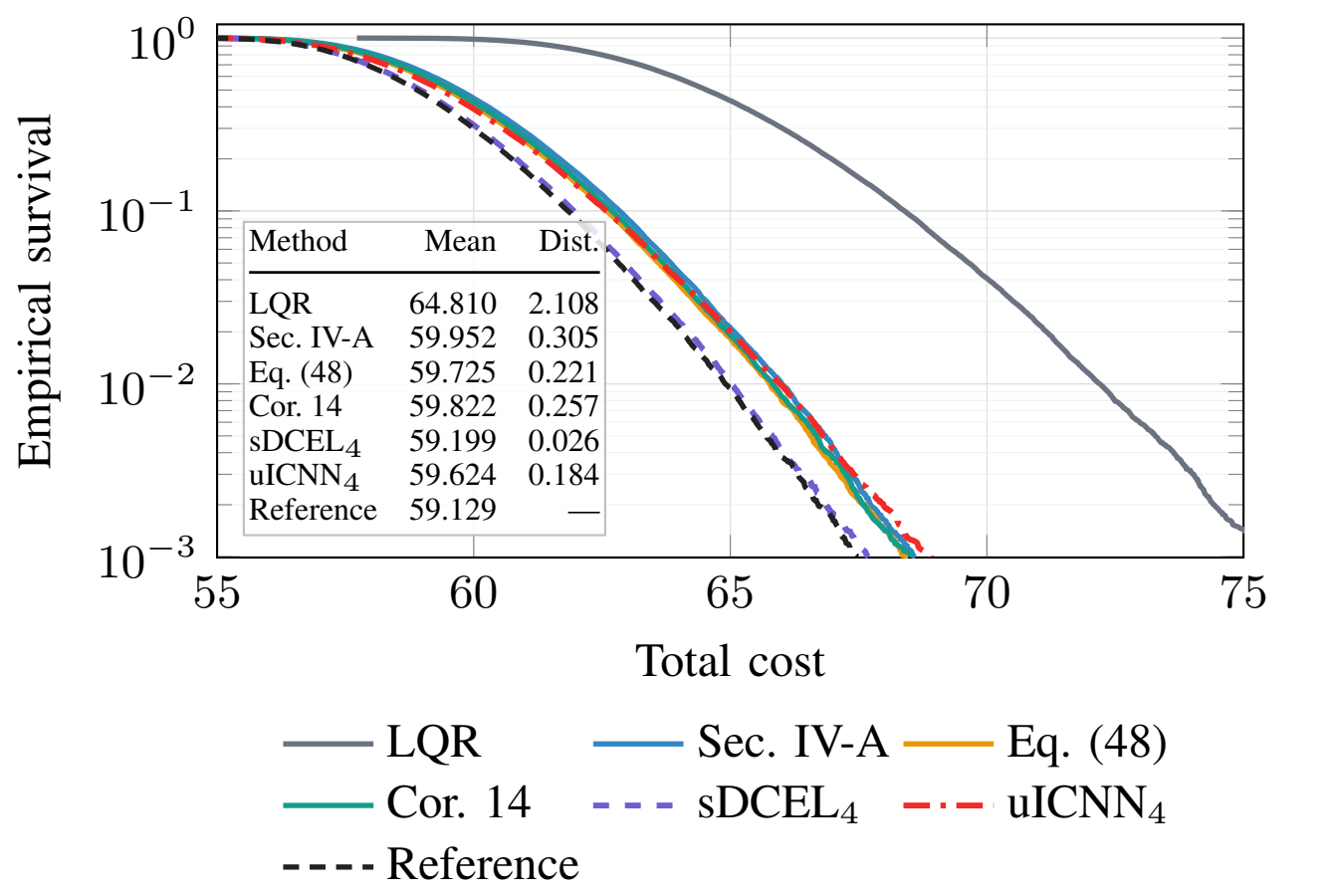

    \paperfigure{adversarial_cost_survival_h15_x2p25_1p5_a09}
    \caption{Empirical total-cost survival functions under each policy's own discrete AVaR adversary, starting from $x_0=(2.25,1.5)^\intercal$ ($j=15$, $\alpha=0.9$, and $50{,}000$ paths per policy). The aDP controllers use $\|x\|_{P_j}^2+r^{\rm dcel}_4(x)$ (sDCEL$_4$) and $g^{\rm icnn}_4(x)$ (uICNN$_4$).}
    \label{fig:adversary-cost-survival}
\end{figure}

\textbf{Evaluation of controllers.}
To evaluate the controllers under adverse disturbance laws,
we simulate \eqref{eq:system-dynamics} with an additive 
disturbance $w_k$ whose law is selected by a (policy-dependent)
adversary. 
For a fixed policy $\bm{\pi}_{j+1}=(\pi_{j+1},\ldots,\pi_1)$,
\eqref{eq:dual_rep_rho_v2} allows to express the risk of the
cost $J_{\bm{\pi}_j}$ as a supremum over disturbance laws.
In the $\AVaR$-Gaussian case we approximate the nominal Gaussian law by
quadrature nodes $w_i$ with probabilities $p_i$, $i=1,\ldots,N_p$.
At a state $x$ with $j+1$ stages remaining, the adversary chooses
\begin{subequations}\label{eq:policy-adversarial-law}
    \begin{align}        
        q_j(x, \bm{\pi}_j)
        {}\in{}&
        \argmax_{\mu \in \mathfrak{M}_\alpha}
        \sum_{i} \mu_i
        J_{\bm{\pi}_j}\bigl(Ax{+}B\pi_{j+1}(x){+}w_i\bigr),
        \\
        \mathfrak{M}_\alpha
        {}={}&
        \left\{
            \mu: \textstyle\sum_i \mu_i=1,\ 0\leq \mu_i\leq\tfrac{p_i}{1-\alpha}
        \right\}.
    \end{align}
\end{subequations}
The stage cost $\ell(x,\pi_{j+1}(x))$ is independent of $w$
and therefore does not affect the maximizer, so it is 
omitted from \eqref{eq:policy-adversarial-law}.
We approximate $J_{\bm{\pi}_j}$ on a state grid, 
draw $w_k$ according to the chosen
probabilities, and repeat the procedure at the next state.
Figures~\ref{fig:adversary-cost-survival-alpha05}
and~\ref{fig:adversary-cost-survival} show the empirical
total-cost survival functions, with the adversarial law
determined separately for each policy.

Figures \ref{fig:adversary-cost-survival-alpha05} and 
\ref{fig:adversary-cost-survival} show that, under each policy's 
own adversarial distribution, the proposed controllers have 
lower empirical mean and upper-tail costs than LQR. 
The gap from LQR is much larger at $\alpha=0.9$: 
the mean of the total cost is $64.810$ for LQR and 
$59.129$ for the reference controller, while the proposed methods 
lie between $59.199$ and $59.952$.
At $\alpha=0.5$, the mean total cost of the proposed methods 
are all close to the reference ($47.027$), compared with 
$48.049$ for LQR. 

To compare each cost distribution with that of the reference method, 
we calculate the empirical Wasserstein-$1$ distance, scaled by the 
interquartile range of the reference distribution \cite{Lu2026}. 
This dimensionless metric measures the average transport discrepancy
as a fraction of the reference distribution's interquartile range.
To define the empirical $W_1$ distance, let $\widehat F_N$ and
$\widehat F_N^{\mathrm{ref}}$ be the empirical \acp{cdf} of two samples
of equal size
$N$, with ordered total costs $C_{(1)}\leq\cdots\leq C_{(N)}$ and
$C^{\mathrm{ref}}_{(1)}\leq\cdots\leq C^{\mathrm{ref}}_{(N)}$, respectively. 
The empirical Wasserstein-$1$ distance is
\begin{equation}
    \widehat W_1
    =\int_{-\infty}^{\infty}|\widehat F_N(t)-\widehat F_N^{\mathrm{ref}}(t)|\mathrm{d}t
    =\tfrac{1}{N}\sum_{i=1}^N
      |C_{(i)}-C^{\mathrm{ref}}_{(i)}|.
\end{equation}

Note that the relative performance of the approximations changes with 
$\alpha$. At $\alpha = 0.5$, the controllers of Section 
\ref{sec:finer-approximation} are slightly closer to the reference 
than the seminorm controller. 
At $\alpha=0.9$, the structured \ac{aDP} controller is closest to the
reference survival curve, including in the upper tail.

\textbf{Performance under common adversary.} It is important to note that each of the methods in 
Figures \ref{fig:adversary-cost-survival-alpha05} and 
\ref{fig:adversary-cost-survival} uses a different, 
policy-dependent adversarial law for $w_k$. 
The figures compare performance against each policy's 
\textit{own} adversary and not a common policy-independent 
adversary.
To evaluate the performance of the difference methods 
against the same adversarial law, we use 
\begin{equation}
    q_j^{\rm com}(x)
    {}\in{}
    \argmax_{q\in\mathfrak{M}_\alpha}
    \sum_i q_i J_{\bm\pi_j^{\rm LQR}}
    \bigl(\bar{A}_{j+1} x + w_i\bigr),
\end{equation}
where $J_{\bm\pi_j^{\rm LQR}}$ is the nested AVaR cost-to-go when following the LQR policy, 
$\pi_{j+1}(x) = K_{j+1} x$.
Roughly speaking, this is a disturbance with a 
distribution chosen adversarially to hurt LQR.

\begin{figure}
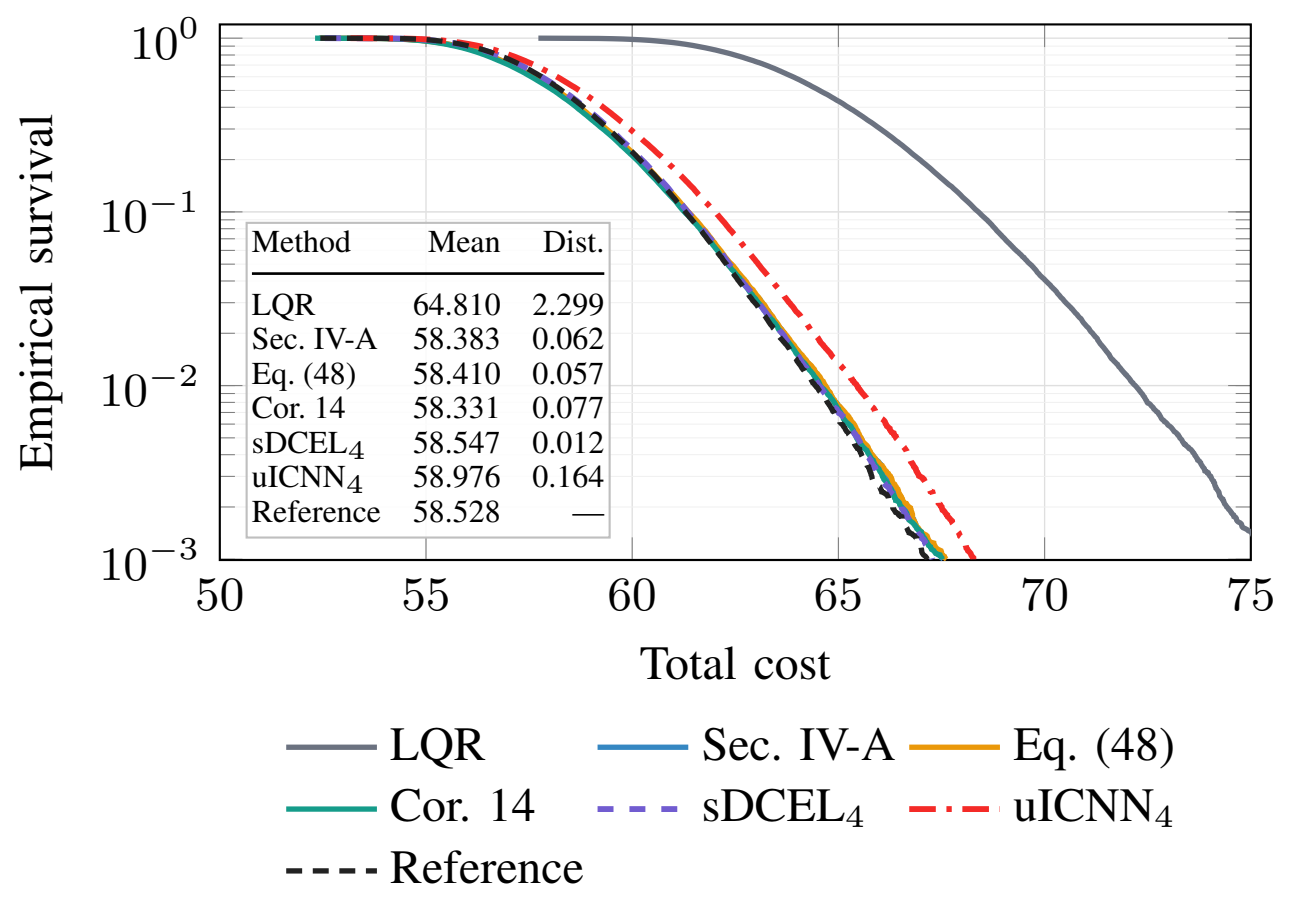

    \paperfigure{adversarial_cost_survival_h15_common_lqr}
    \caption{Total-cost survival under a common LQR-based adversary, starting from $x_0=(2.25,1.5)^\intercal$ ($j=15$, $\alpha=0.9$, and $50{,}000$ paths per policy).}
    \label{fig:common-lqr-adversary-cost-survival}
\end{figure}

Figure~\ref{fig:common-lqr-adversary-cost-survival} shows that, under this common adversary, the proposed controllers have substantially lower mean costs than LQR: their means range from $58.33$ to $58.55$, compared with $64.81$ for LQR. Their upper tails are also much lighter, and their survival curves lie close to that of the reference controller.

Under this common adversary, structured aDP reproduces 
the reference total cost distribution more closely than unstructured aDP, 
with IQR-normalized $W_1$-distances of $0.012$ and $0.164$, respectively.

% SDDP
% ------------
% Key message: SDDP is more accurate, as expected,
% but significantly slower compared to Prop 7.

% Note that the popular SDDP numerical 
% method \cite{FreitasDaCosta2023,GuiguesRomisch2012} requires a discretization of the continuous 
% distribution of $w$.

\subsection{System with exponential disturbance}
Consider the same dynamical system as in the 
previous section, but with a heavier-tailed 
disturbance with \ac{pdf} $p_w(w) \propto \exp(-\|w\|_S)$,
with $S = \left[\begin{smallmatrix}
   300 & -150 \\ -150 & 450
\end{smallmatrix}\right]$.
Table~\ref{tab:value-relative-error-exponential} reports the
relative errors obtained using the same grid and error metric 
as in the preceding subsection.
The reference method is the same as in Section 
\ref{sec:simulations:gaussian_dist} but with a 
generalized Gauss-Laguerre quadrature
\cite[Sec.~3.5(v)]{Olver2010}.
The entropic value-at-risk does not yield a well-defined 
problem (except for trivial cases such as $\alpha=0$):
the distribution of $w$ is such that $\exp(t\|Ax+Bu+w\|_{P_f}^2)$ 
has non-finite expectation for all $t>0$ and for any $x$ and $u$
and $P_f\succ 0$, so $V_j^\star(x)$ become infinite for all $j\geq 1$.

\begin{SCtable}[1][htpb]
    \scriptsize
    \setlength{\tabcolsep}{4.2pt}
    \renewcommand{\arraystretch}{1.08}
\begin{tabular}{cccc@{\hspace{10pt}}c}
\toprule
$\alpha$ & $j$ & Sec.~\ref{sec:sum-seminorms} & Eq.~\eqref{eq:ab_tail_expectations} & LQ \\
\midrule
0.10 & 1  & \heatA{0.1} & \heatA{0.0} & \heatD{1.4} \\
     & 5  & \heatC{0.6} & \heatA{0.1} & \heatE{3.8} \\
     & 10 & \heatC{0.8} & \heatA{0.1} & \heatE{4.0} \\
     & 20 & \heatD{1.1} & \heatB{0.2} & \heatE{4.1} \\
\addlinespace[0.3ex]
0.50 & 1  & \heatB{0.3} & \heatA{0.0} & \heatE{5.2} \\
     & 5  & \heatD{1.5} & \heatC{0.5} & \heatF{12.8} \\
     & 10 & \heatD{2.6} & \heatC{0.9} & \heatF{13.3} \\
     & 20 & \heatE{3.5} & \heatD{1.1} & \heatF{13.4} \\
\addlinespace[0.3ex]
0.90 & 1  & \heatC{0.6} & \heatB{0.2} & \heatF{11.9} \\
     & 5  & \heatE{3.9} & \heatD{2.3} & \heatF{26.1} \\
     & 10 & \heatE{7.0} & \heatE{3.5} & \heatF{26.4} \\
     & 20 & \heatE{9.1} & \heatE{4.2} & \heatF{26.6} \\
\addlinespace[0.3ex]
0.95 & 1  & \heatC{0.8} & \heatB{0.3} & \heatF{14.3} \\
     & 5  & \heatE{5.2} & \heatE{3.1} & \heatF{30.2} \\
     & 10 & \heatE{9.1} & \heatE{4.7} & \heatF{30.3} \\
     & 20 & \heatF{11.5} & \heatE{5.6} & \heatF{30.5} \\
\bottomrule
\end{tabular}
% \par\vspace{0.5ex}\hspace{3em}
% {\small lower error\enspace
% \heatlegendbox{HeatA}\heatlegendbox{HeatB}\heatlegendbox{HeatC}%
% \heatlegendbox{HeatD}\heatlegendbox{HeatE}\heatlegendbox{HeatF}%
% \enspace higher error}

    \caption{Relative value-function error as defined in \eqref{eq:value-relative-error} for the exponential
    disturbance, for different values of the $\AVaR$ level $\alpha$ and
    horizon length $j$.}
    \label{tab:value-relative-error-exponential}
\end{SCtable}

\subsection{Higher-dimensional system}%
\label{sec:ten-state-system}
We consider a dynamical system of the form 
\eqref{eq:system-dynamics} with $n_x=10$ states 
and $n_u=5$ inputs. The system comprises five weakly 
coupled subsystems; the system matrices are
\begin{equation}
    A = 
    \left[
    \begin{smallmatrix}
        A_0 & D
        \\ 
        & A_0 & D
        \\ 
        & & A_0 & D
        \\ 
        & & & A_0 & D
        \\ 
        D & & & & A_0 & D
    \end{smallmatrix}
    \right],
    B = I_5 \otimes B_0,
\end{equation}
where 
\begin{equation}
    A_0=\left[\begin{smallmatrix}0.88&0.18\\-0.04&0.82\end{smallmatrix}\right],
    B_0=\left[\begin{smallmatrix}0.12\\0.42\end{smallmatrix}\right],
    D = \left[\begin{smallmatrix}0.025 & 0\\0 & 0\end{smallmatrix}\right].
\end{equation}
The weight matrices are $Q_0=I_5\otimes\diag(1,0.35)$,
$R=0.25I_5$ and choose $P_f$ as the solution of the
discrete-time algebraic Riccati equation.
The disturbances, $w_k$, are iid zero-mean Gaussian with 
covariance matrix $\Sigma_{rs}=0.018 \cdot 0.35^{|r-s|}$, 
for $r,s=1,\ldots,10$.
To obtain a system with dense matrices, we consider the 
change of variables $x'=Ux$, and $w'=Uw$, where 
$U \in \R^{10\times 10}$ is a random orthogonal matrix.

\begin{figure}[htpb]
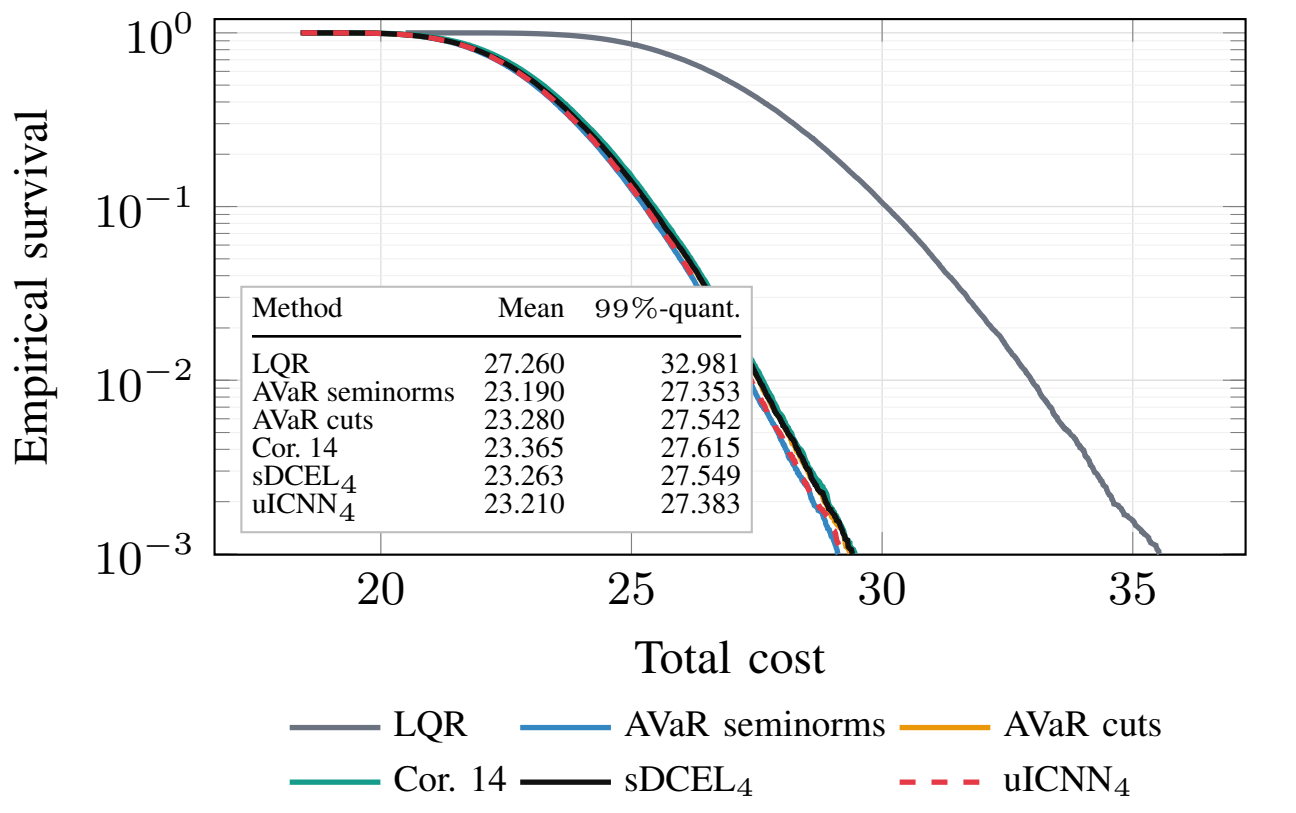

    \centering
    \paperfigure{adversarial_cost_survival_nx10_h10_energy3_a09}
    \caption{Total-cost survival for the ten-state system under a common LQR-guided Gaussian-tail stress law ($N=10$, $\alpha=0.9$, $\|x_0\|_{P_f}=3$; $50{,}000$ paths per controller). Both neural controllers use $\AVaR_{0.9}$.}
    \label{fig:ten-state-common-lqr-cost-survival}
\end{figure}

\textbf{Value function approximations.}
For the ten-state system, computing a reference optimal
value function is computationally prohibitive due to
the dimensionality of the state space.
We therefore assess the closed-loop performance of the
proposed controllers through the empirical total-cost
survival functions shown in
Fig.~\ref{fig:ten-state-common-lqr-cost-survival}.
All proposed methods substantially reduce the
probability of high-cost outcomes compared with the
risk-neutral LQR controller, demonstrating their
risk-averse behavior.
Moreover, the proposed methods exhibit similar
performance, with the sum-of-seminorms approximation
showing a slight advantage.

\begin{figure}[htpb]
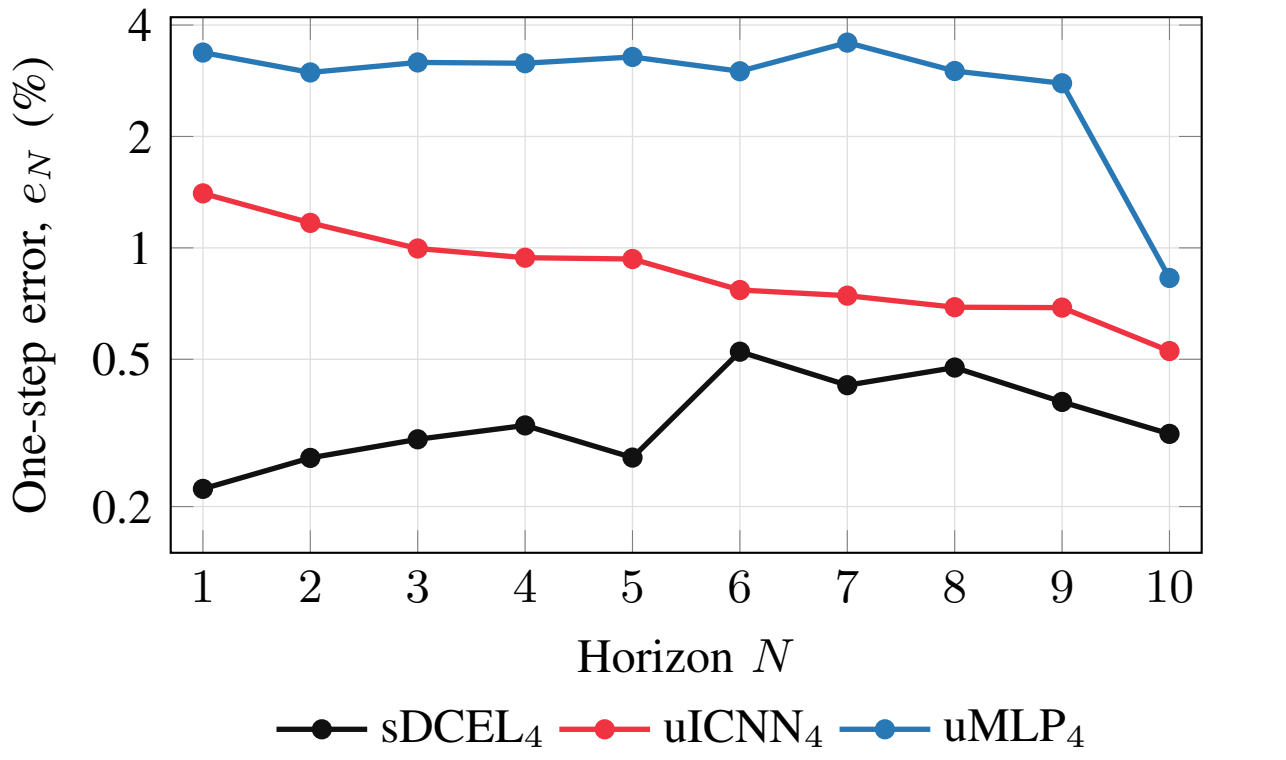

    \centering
    \paperfigure{bellman_rms_scaled_nx10_energy6p25}
    \caption{Bellman errors $e_N$ for the system 
    of Section \ref{sec:ten-state-system} with $\|x\|_{P_f}=6.25$, $\alpha=0.9$, evaluated on $2{,}048$ test points.}
    \label{fig:ten-state-bellman-rms-extrapolation}
\end{figure}

\textbf{Approximate DP.}
To assess the aDP approximations on unseen states without computing
$V_N^\star$, we measure consistency with the Bellman recursion.
For held-out test states $x_1,\ldots,x_{n_{\mathrm{test}}}$, let
$y_{N,i}\approx(\T\widetilde V_{N-1})(x_i)$ denote the numerical
Bellman targets. We define the relative Bellman error as
\begin{equation}
    e_N
    \coloneqq
    \sqrt{
        \frac{\sum_{i=1}^{n_{\mathrm{test}}}
        |\widetilde V_N(x_i)-y_{N,i}|^2}
        {\sum_{i=1}^{n_{\mathrm{test}}}|y_{N,i}|^2}
    } {}\cdot{} 100\%.
    \label{eq:rms-scaled-bellman-error}
\end{equation}
This measures the one-step Bellman residual relative to the magnitude
of the target costs and is reported in
Fig.~\ref{fig:ten-state-bellman-rms-extrapolation}.
This is a consistency metric rather than a
measure of the error relative to $V_N^\star$.

\textbf{Computational complexity.}
The sum-of-seminorms solves online an SOCP problem 
with 16 decision variables and 11 second-order 
conic constraints. MOSEK~\cite{MOSEK1127} yields runtimes
with $95\%$-quantile $\unit[1.8]{ms}$ based on 
384 solved problems.
The supporting hyperplanes method of 
\eqref{eq:ab_tail_expectations} and 
Proposition \ref{prop:elliptic-symmetry-cut}
as well as Corollary \ref{cor:gaussian-nontight}
solve an SOCP with 116 variables and 1,490 
linear inequalities; the $95\%$-quantile is 
$\unit[8.8]{ms}$.
Regarding the aDP methods, the $95\%$-quantiles 
of the runtimes of sDCEL$_4$, uICNN$_4$, and 
uMLP$_4$ are $\unit[11.1]{ms}$,
$\unit[22.2]{ms}$, and 
$\unit[53.6]{ms}$, respectively.

% {\color{gray}\textbf{To be removed:}
% We briefly describe the network classes used below. 
% A standard $L$-layer \ac{MLP} is written as%
% \begin{subequations}%
%     \label{eq:adp-mlp}
%     \begin{align}
%         h^{(0)}(x)&=x,\\
%         h^{(k+1)}(x)&=\sigma\bigl(W_kh^{(k)}(x)+b_k\bigr),\\
%         F_L(x)&=w_L^\intercal h^{(L)}(x)+b_L,
%     \end{align}
% \end{subequations}
% where $\sigma(t)=t(1+\exp(-t))^{-1}$ is the swish function.
% An even \ac{MLP} can be constructed by symmetrization 
% using $r^{\rm even}_L(x) = \tfrac12\bigl(F_L(x)+F_L(-x)\bigr)$.

% For $\tau>0$, define the following log-sum-exp-type function
% \begin{equation*}
%     \chi_{\tau}(z_1,\ldots,z_p)
%     =\tau\log \left(\tfrac{1}{p}\sum_{i=1}^{p}\exp(z_i/\tau)\right).
% \end{equation*}
% Set $s_\ell^{\pm}(x)=b_\ell\pm a_\ell^\intercal x$. For the log-sum-exp
% and the max-affine residuals, respectively, we define
% \begin{align}
%     r^{\rm lse}(x)
%     &=\chi_{\tau}\!\bigl(
%         c,\bigl(s_\ell^+(x),s_\ell^-(x)\bigr)_{\ell=1}^{M}
%         \bigr),
%     \label{eq:adp-lse}
%     \\
%     r^{\rm ma}(x)
%     &=\max\bigl\{c,
%         \max_{1\leq\ell\leq M}
%         \bigl\{b_\ell+|a_\ell^\intercal x|\bigr\}\bigr\}.
%     \label{eq:adp-max-affine}
% \end{align}
% The latter is a symmetric maxout-type network \cite{Goodfellow+2013}.
% }

\section{Conclusions and future work}
We studied finite-horizon risk-averse \ac{LQ} optimal control
under nested law-invariant coherent risk measures.
Our main result establishes that, although the optimal
value functions are generally nonquadratic, they retain
the classical Riccati-generated quadratic component,
augmented by a convex, even, and globally Lipschitz residual.
This decomposition admits a regularized \ac{LQ} interpretation,
where the residual captures the additional cost induced
by risk aversion while preserving the underlying Riccati structure.
Furthermore, the dual representation of coherent risk measures
provides a game-theoretic interpretation of the multistage
risk-averse problem, in which the controller acts against
an adversarial selection of probability distributions.
Exploiting the identified structure, we developed tractable
value-function approximations with global error bounds
and finite-horizon policy suboptimality guarantees.
We further improved these approximations using
supporting hyperplanes, with explicit constructions
for $\AVaR$ and $\EVaR$.
Numerical results demonstrated that these approximations
achieve high accuracy with low computational effort.
Moreover, the proposed structure-preserving \ac{aDP}
approach outperformed unstructured alternatives in
value-function approximation, highlighting the benefits
of incorporating analytical structure into numerical
approximation methods.

Future work will develop recursive supporting cuts for $\rho[\xi_j(m,w)+\widehat{r}_j(m+w)]$, 
thereby accounting explicitly for the disturbance-dependent residual. 
This will connect the decomposition of 
Theorem \ref{thm:Vjstar_structure} 
to risk-averse SDDP \cite{GuiguesRomisch2012}. 
Establishing global error bounds, 
guaranteeing \hyperlink{D1}{D1} and \hyperlink{D2}{D2}, and 
assessing the computational cost
are open questions to be addressed.

\bibliographystyle{IEEEtran}
\bibliography{bibliography}

\appendix
\hypertarget{link:appendix}{}
\begin{proof}[Proof of Lemma \ref{lem:psi_properties}]
    (\textit{Finiteness}) Let $\zeta(m,w)=2m^\intercal Pw+\|w\|_P^2+r(m+w)$. 
    % 1. ρ(|w|) is finite
    Since $\|w\| \leq \tfrac{1}{2}(\|w\|^2 + 1)$,
    by applying $\rho$ to both sides and by 
    \hyperlink{A1}{A1}, \hyperlink{A2}{A2}, 
    and \hyperlink{A3}{A3},
    $\rho(\|w\|) \leq \tfrac{1}{2}(\rho(\|w\|^2) + 1) < \infty$, by Assumption \ref{assm:rho_disturbance}.
    % 2. The terms in ζ are bounded by ||w|| and ||w||^2
    Since $r$ is $L_r$-Lipschitz,
    $|r(m+w)|\leq |r(m)| + L_r\|w\|$.
    By the Cauchy-Schwarz, inequality
    $|m^\intercal Pw| \leq \|Pm\|\,\|w\|$.
    Moreover, $\|w\|_P^2 \leq \|P\|\,\|w\|^2$,
    so $|\zeta(m,w)| \leq  (2\|Pm\| + L_r)\|w\| + \|P\|\,\|w\|^2+|r(m)|$.
    % 3. ψ is finite
    By the coherence axioms of $\rho$,
    $\psi(m) =  \rho(\zeta(m,w)) \leq \rho(|\zeta(m,w)|) \leq (2\|Pm\|{+}L_r)\rho(\|w\|) {+}\|P\|\rho(\|w\|^2) {+} |r(m)| < \infty$.

    (\textit{Convexity}) For fixed $w$, $\zeta({}\cdot{},w)$ is convex, therefore monotonicity and convexity of $\rho$ imply
    that, for $m_1, m_2\in\R^{n_x}$ and $\theta\in [0,1]$,
    and using the notation $\bar{m} = \theta m_1 + (1-\theta)m_2$,
    $\zeta_1 = \zeta(m_1, w)$, $\zeta_2 = \zeta(m_2, w)$,
    \begin{subequations}    
    \begin{align}
        &\psi(\bar{m}) 
        {}={}
        \rho(\zeta(\bar{m}, w))
        \\
        &\quad{}\leq{}
        \rho(\theta\zeta_1 + (1-\theta)\zeta_2)
        \label{eq:psi:2}
        \\
        &\quad{}\leq{}
        \rho(\theta\zeta_1)+ \rho((1-\theta)\zeta_2)
        \label{eq:psi:3}
        \\
        &\quad{}={}
        \theta\rho(\zeta_1)+ (1-\theta)\rho(\zeta_2)
        \label{eq:psi:4}
        \\
        &\quad{}={}
        \theta\psi(m_1)+ (1-\theta)\psi(m_2),
        \label{eq:psi:5}
    \end{align}
    \end{subequations}
    where \eqref{eq:psi:2} is due to the convexity 
    of $\zeta({}\cdot{}, w)$ and monotonicity of $\rho$, 
    \eqref{eq:psi:3}
    is by the subadditivity of $\rho$, 
    and \eqref{eq:psi:4} is by positive homogeneity
    of $\rho$.

    (\textit{Evenness})
    Since $w\overset{\mathrm d}{\sim}-w$, law invariance of
    $\rho$ \cite[Sect. 6.3.3]{ShapiroDentcheva2021} and
    the evenness of $r$ give
    \begin{align*}
        \psi(-m)
        &=\rho\left(\zeta(-m,-w) \right) \\
        &=\rho\left(2m^\intercal Pw+\|w\|_P^2+r(m+w)\right)
          =\psi(m).
    \end{align*}

    (\textit{Lipschitzianity})
    Let $X, Y \in \mathcal{Z}$ with $\rho(|X|), \rho(|Y|) < \infty$; then from \hyperlink{A1}{A1} and 
    \hyperlink{A4}{A4}, $\rho(X)-\rho(Y) \leq \rho(X-Y) \leq \rho(|X-Y|)$. Likewise, 
    $\rho(Y)-\rho(X) \leq \rho(Y-X) \leq \rho(|X-Y|)$,
    therefore 
    \begin{equation}
        |\rho(X)-\rho(Y)|
        \leq
        \rho(|X-Y|).
        \label{eq:rho_useful_ineq}
    \end{equation}
    Let $m,m'\in\R^{n_x}$; we have
    \begin{align}
        &|\zeta(m,w){-}\zeta(m',w)|
        \notag
        \\
        &\quad{}={}
        \left|
            2(m{-}m')^\intercal Pw
            {+}r(m{+}w){-}r(m'{+}w)
        \right|
        \notag
        \\
        &\quad{}\leq{}
        \bigl(2\|Pw\|+L_r\bigr)\|m-m'\|.
    \end{align}
    Applying \eqref{eq:rho_useful_ineq} with
    $X=\zeta(m,w)$ and $Y=\zeta(m',w)$, and using monotonicity,
    positive homogeneity, and translation equivariance of $\rho$, gives
    \begin{align*}
        |\psi(m)-\psi(m')|
        {}\leq{}&
        \rho\bigl(|\zeta(m,w)-\zeta(m',w)|\bigr)\\
        {}\leq{}&
        \rho\left(
            \|m-m'\|\bigl(2\|Pw\|+L_r\bigr)
        \right)\\
        {}={}&
        \left(
            2\rho(\|Pw\|)+L_r
        \right)
        \|m-m'\|,
    \end{align*}
    where
    $L_\psi \coloneqq 2\rho(\|Pw\|)+L_r <\infty$,
    hence $\psi$ is globally Lipschitz.
\end{proof}

\begin{proof}[Proof of Prop. \ref{prop:quadratic-centered-error}]
    \hypertarget{link:appendix:proof-quadratic-upper-bound}{}
    From Theorem \ref{thm:Vjstar_structure}, 
    $V_j^\star(x) = \|x\|_{P_j}^2 + r_j(x)$, where $r_j$ is 
    convex, even, and $L_{r_j}$-Lipschitz. Let 
    $\Theta_j(m) {}\coloneqq{} \rho(V_j^\star(m + w))$.

    Firstly, let us briefly describe the sign trick of 
    \cite[Lem. 6.6.2]{Vershynin2026}: Let $g\sim \mathcal{N}(0,1)$. We can write $g = \sign(g) |g|$; 
    let $\varepsilon {}\coloneqq{} \sign(g)$. Then, 
    by the independence of $\varepsilon$ and $|g|$,
    $\E[g {}\mid{} \varepsilon] = \varepsilon \E[|g|] = \varepsilon \sqrt{\nicefrac{2}{\pi}}$.
    % , so 
    % $\varepsilon = \E[\sqrt{\nicefrac{\pi}{2}} g {}\mid{} \epsilon]$.
    
    Without loss of generality 
    suppose $m\neq0$ and let $w'\sim\mathcal{N}(0, \Sigma_w)$,
    independent of $w$. We apply the sign trick to
    \begin{equation}
        g = \frac{m^\intercal \Sigma_w^{-1} w'}{\|\Sigma_w^{-1/2}m\|},
    \end{equation}
    with $\varepsilon = \sign(g) = \sign(m^\intercal \Sigma_w^{-1} w')$.
    
    We then have 
    \begin{subequations}        
        \begin{align}
            &\E[w' {}\mid{} \varepsilon] 
            {}={}
            \varepsilon \sqrt{\tfrac{2}{\pi}} \frac{m}{\|\Sigma_w^{-1/2}m\|}
            \\
            \Rightarrow{}
            &
            \E\left[w + \sqrt{\tfrac{\pi}{2}}\|\Sigma_w^{-1/2}m\|w' {}\;\Big|\;{} w, \varepsilon\right]
            {}={} w + \varepsilon m.
        \end{align}
    \end{subequations}
    Let us define 
    $Y = w + \sqrt{\tfrac{\pi}{2}}\|\Sigma_w^{-1/2}m\|w'$.
    Since $V_j^\star$ is convex, the conditional Jensen's 
    inequality \cite[p.~143]{Cinlar2011} gives\footnote{%
            The integrability conditions of the conditional 
            Jensen inequality are satisfied: indeed, 
            from Theorem~\ref{thm:Vjstar_structure},
    $|V_j^\star(x)| \leq \|P_j\|_2\|x\|^2 + L_{r_j}\|x\| + |r_j(0)|$,
    and $\E\|Y\|^2 < \infty$ since $Y$ is Gaussian, 
    so $\E|V_j^\star(Y)| < \infty$.. 
        } 
    \begin{equation}
        V_j^\star(
                \underbracket[0.5pt]{
                    \E[Y \mid w, \varepsilon]
                }_{w + \varepsilon m}
        )
        {}\leq{}
        \E[V_j^\star(Y) {}\mid{} w, \varepsilon],
    \end{equation}
    a.s., and by virtue of \hyperlink{A1}{A1},
    \begin{equation}
        \rho[V_j^\star(w + \varepsilon m)]
        {}\leq{}
        \rho[\E[V_j^\star(Y) {}\mid{} w, \varepsilon]]
        {}\overset{\text{DM}}{\leq}{}
        \rho[V_j^\star(Y)].
        \label{eq:Theta-upper-bound-1}
    \end{equation}
    Since $w \eqd -w$ and $V_j^\star$ is even,
    $V_j^\star(w-m) \eqd V_j^\star(-w-m) = V_j^\star(w+m)$, so 
    $V_j^\star(w+\varepsilon m) \eqd V_j^\star(w+m)$, since $\varepsilon$ is independent of $w$. By the law invariance of $\rho$, 
    $\Theta_j(m) = \rho[V_j^\star(w + \varepsilon m)]$.
    As a result, \eqref{eq:Theta-upper-bound-1} becomes 
    $\Theta_j(m) \leq \rho(V_j^\star(Y))$.

    By the definition of $Y$ we have 
    $Y \sim \mathcal{N}(0, a^2 \Sigma_w)$, 
    where $a^2 = 1 + \tfrac{\pi}{2}\|\Sigma_w^{-1/2}m\|^2$.
    Equivalently, $Y \eqd aw$, so by the law invariance
    property, $\Theta_j(m) \leq \rho(V_j^\star(aw))$.
    Exploiting the Lipschitzianity of $r_j$, 
    \begin{align}
        V_j^\star(aw) - V_j^\star(w)
        {}={}&
        (a^2-1)\|w\|_{P_j}^2 + r_j(aw) - r_j(w)
        \notag\\
        {}\leq{}&
        (a^2-1)\|w\|_{P_j}^2 + L_{r_j}(a-1)\|w\|.
        \label{eq:Vstar-quad-upper-bound}
    \end{align}
    By applying $\rho$ to both sides of \eqref{eq:Vstar-quad-upper-bound}
    and by virtue of \hyperlink{A1}{A1}, \hyperlink{A3}{A3}, 
    and \hyperlink{A4}{A4}
    \begin{align}
        \Theta_j(m)-\Theta_j(0)
        {}\leq{}&
        (a^2-1)\rho(\|w\|_{P_j}^2)
        +L_{r_j}(a-1)\rho(\|w\|)
        \notag
        \\
        {}\leq{}&
        \tfrac{\pi}{2}\|\Sigma_w^{-1/2}m\|^2\rho(\|w\|_{P_j}^2)
        \notag
        \\
        &\quad{}+{}
        L_{r_j}\tfrac{\pi}{4}\|\Sigma_w^{-1/2}m\|^2 \rho(\|w\|),
    \end{align}
    where we used the fact that $a\geq 1$, 
    thus 
    \begin{equation}
        a-1=\frac{a^2-1}{a+1}\leq\frac{a^2-1}{2}\leq\frac{\pi}{4}\|\Sigma_w^{-1/2}m\|^2.
    \end{equation}
    Let us define the constant 
    $b_j {}\coloneqq{} \tfrac{\pi}{2}
        \bigl[ 
        \rho(\|w\|_{P_j}^2)
        +\tfrac{1}{2}L_{r_j}\rho(\|w\|)
        \bigr]$.
    Since $\Theta_j$ is convex and even, $\Theta_j(m) \geq \Theta_j(0)$,
    so
    \begin{equation}
        0
        {}\leq{}
        \Theta_j(m) - \Theta_j(0)
        {}\leq{}
        b_j
        \|\Sigma_w^{-1/2}m\|^2.    
    \end{equation}
    Next, note that $\psi_j(m) {}={} \Theta_j(m)-\|m\|_{P_j}^2$.

    Using $v=0$ in \eqref{eq:r-j-plus-one} and since 
    $r_{j+1}(0) = \psi_j(0) = \Theta_j(0)$, we have 
    \begin{align}
        0 
        {}\leq{}
        r_{j+1}(x)-r_{j+1}(0)
        {}\leq{}&
        \psi_j(\bar{A}_{j+1}x) - \psi_j(0)
        \notag
        \\
        {}\leq{}&
        \Theta_j(\bar{A}_{j+1}x) - \Theta_j(0)
        \notag
        \\
        {}\leq{}& 
        b_j \|\Sigma_w^{-1/2}\bar{A}_{j+1}x\|^2.
    \end{align}
    By the definition of $\eta_j$ in \eqref{eq:eta_j:definition},
    \begin{align}
        \eta_j(m) - \eta_j(0)
        {}={}&
        \rho(\|m+w\|_{P_j}^2) - \rho(\|w\|_{P_j}^2) 
            - \|m\|_{P_j}^2
        \notag
        \\
        {}\leq{}&
        \rho(\|m+w\|_{P_j}^2) - \rho(\|w\|_{P_j}^2)
        \notag
        \\
        {}\leq{}&
        \tfrac{\pi}{2}\rho(\|w\|_{P_j}^2)
        \|\Sigma_w^{-1/2}m\|^2,
        \label{eq:delta-eta-quad-upper-bnd}
    \end{align}
    where the last inequality follows by applying the 
    above argument to the function $z\mapsto\|z\|_{P_j}^2$, 
    whose residual is zero, so the upper bound 
    in \eqref{eq:delta-eta-quad-upper-bnd} follows by using 
    $L_{r_j}=0$ in the definition of $b_j$.
    From Proposition \ref{prop:eta-cut-approximation} and the
    above bound,
    \begin{equation}
        0
        {}\leq{}
        \underline{\eta}_j(m)-\underline{\eta}_j(0)
        {}\leq{}
        \eta_j(m)-\eta_j(0)
        {}\leq{} 
        b_j\|\Sigma_w^{-1/2}m\|^2.
    \end{equation}
    Now let us define $c_0=0$ and, recursively,
    \begin{equation}
        c_{j+1}
        =c_j\|\bar{A}_{j+1}\|_2^2
        +b_j\|\Sigma_w^{-1/2} \bar{A}_{j+1}\|_2^2.
    \end{equation}
    Since $s_0^*=0$, induction using \eqref{eq:s-cut} gives
    $0 {}\leq{} s_j^*(x)-s_j^*(0) {}\leq{} c_j\|x\|^2$.
    Indeed, if this holds at stage $j$, then
    \begin{align}
        0
        {}\leq{}& 
        s_{j+1}^*(x) - s_{j+1}^*(0)
        \notag
        \\
        {}\leq{}&
        c_j\|\bar{A}_{j+1}x\|^2
        + b_j \|\Sigma_w^{-1/2}\bar{A}_{j+1}x\|^2
        \leq c_{j+1}\|x\|^2.
    \end{align}
    The bound already established for $r_{j+1}$ also gives
    $0 {}\leq{} r_{j+1}(x) - r_{j+1}(0) {}\leq{} c_{j+1}\|x\|^2$,
    so
    \begin{align}
        |e_{j+1}(x)|
        {}={}&
        \bigl|
            s_{j+1}^*(x)
            - s_{j+1}^*(0)
            - r_{j+1}(x)
            + r_{j+1}(0)
        \bigr|
        \notag 
        \\
        {}\leq{}& 
        c_{j+1}\|x\|^2.
    \end{align}
    This completes the proof.
\end{proof}

% This is only used when `make figures` is run
\ifdefined\ExtraPaperFigures
    \clearpage\ExtraPaperFigures
\fi
\end{document}